\documentstyle{amsppt}  
\magnification1200
\pagewidth{6.5 true in}
\pageheight{9.25 true in}
\NoBlackBoxes

\def\phi{\varphi}

\def\L{\fracwithdelims()}
\def\sumf{\sideset \and^\flat \to \sum} 
\def\sumstar{\sideset \and^* \to \sum}
\def\lam{\lambda}
\topmatter
\title The distribution of values of $L(1,\chi_d)$
\endtitle
\author Andrew Granville and K. Soundararajan
\endauthor
\address{Department of Mathematics, University of Georgia, Athens, 
Georgia 30602, USA} \endaddress
\email{andrew{\@}math.uga.edu} \endemail
\address{Department of Mathematics, University of Michigan, Ann Arbor, 
Michigan 48109, USA} \endaddress
\email{ksound{\@}umich.edu} \endemail
\thanks{Both authors are partially supported by the 
National Science Foundation.  The second 
author is partially supported by the American Institute of Mathematics (AIM).}
\endthanks
\endtopmatter

\head  1. Introduction \endhead 

\noindent Throughout this paper $d$ will denote a fundamental discriminant, 
and $\chi_d$ the associated primitive real character to the modulus $|d|$.  
We investigate here the distribution of values of $L(1,\chi_d)$ 
as $d$ varies over all fundamental discriminants with $|d|\le x$.  
Our main concern is to compare the distribution of values of $L(1,\chi_d)$ 
with the distribution of ``random Euler products" $L(1,X) = \prod_{p} 
(1-X(p)/p)^{-1}$ where the $X(p)$'s are independent random variables 
taking values $0$ or $\pm 1$ with suitable probabilities, described below.
For example, we shall give asymptotics for the probability that 
$L(1,\chi_d)$ exceeds $e^{\gamma} \tau$, and the probability that 
$L(1,\chi_d) \le \frac{\pi^2}{6}\frac{1}{e^{\gamma} \tau}$ uniformly in 
a wide range of $\tau$.  These results are sufficiently uniform to
prove slightly more than a recent conjecture of H.L.~Montgomery and 
R.C.~Vaughan [15].  One important motivation for our work is to make 
progress towards resolving the discrepancy between extreme 
values that may be exhibited (the omega results of 
S.D. Chowla described below) and the conditional bounds on these 
extreme values (the $O$-results of J.E. Littlewood, see below).  
The uniformity of our results provides evidence that the 
omega results of Chowla may represent the true nature of extreme values of
$L(1,\chi_d)$.

These questions have been studied by many authors, most notably by P.D.T.A. 
Elliott [5,6,7] and Montgomery and Vaughan [15].  We begin by reviewing 
some of the history of the subject which will help place our results in 
context.  Throughout this paper $\log_j$ will denote the $j$-fold iterated 
logarithm; that is, $\log_2 =\log \log$, $\log_3 =\log \log \log $, and so 
on.  In [11] Littlewood showed that on the assumption of 
the Generalized Riemann Hypothesis (GRH) 
$$
\Big(\frac 12 +o(1)\Big) \frac{\zeta(2)}{e^{\gamma} \log_2 |d|} 
\le L(1,\chi_d) \le (2+o(1)) e^{\gamma} \log_2 |d|. \tag{1.1}
$$
He also showed (again assuming GRH) that there are infinitely many 
$d$ such that $L(1,\chi_d) \ge (1+o(1)) e^{\gamma} \log_2 |d|$, and 
that for infinitely many $d$, $L(1,\chi_d)\le (1+o(1)) \zeta(2)/(e^{\gamma} 
\log_2|d|)$.  The latter result was later established unconditionally 
by Chowla [2].   Thus only a factor of $2$ remains at issue regarding the 
extreme values of $L(1,\chi_d)$, under the assumption of GRH.  
This question was addressed in a numerical study 
of D. Shanks [16], but the data there is inconclusive.  Recently 
R.C. Vaughan [17] and H.L. Montgomery and 
Vaughan [15] have returned to this problem, and initiated a 
finer study of these extreme values.  

Write 
$$
\log L(1,\chi\sb d)
=\sum\sb p{\chi\sb d(p)\over p}
+\sum\sb p\sum\sp \infty\sb {k=2}{\chi\sb d(p)\sp k\over kp\sp k}.
$$
The second sum above is rapidly convergent, and hence easy to understand. 
For a typical $d$ one may expect that the first sum 
above behaves like $\sum\sb pX\sb p/p$, where the $X\sb p$ are independent
random variables taking values $\pm 1$ with equal probability. Pursuing 
this probabilistic model, Montgomery and Vaughan (developing ideas
of Montgomery and Odlyzko [14]) suggest that the 
proportion of fundamental discriminants $|d|\le x$ with 
$L(1,\chi_d) > e^{\gamma} \tau$ (say) lies between $\exp(-C e^{\tau}/\tau)$, 
and $\exp(-c e^{\tau}/\tau)$ for appropriate constants $0 < c <C<\infty$.  
A similar assertion holds for the frequency with which 
$L(1,\chi_d) < \zeta(2) e^{-\gamma}/\tau$.
Extrapolating this model, they formulated three conjectures 
on the frequencies with which certain extreme values occur.  

\proclaim{Conjecture 1 (Montgomery-Vaughan)} 
The proportion of fundamental discriminants 
$|d|\le x$ with $L(1,\chi_d) \ge e^{\gamma} \log_2 |d|$ 
is $> \exp(-C\log x/\log_2 x)$ and $< \exp(-c\log x/\log_2 x)$ 
for appropriate constants $0<c<C<\infty$.  Similar estimates 
apply to the proportion of fundamental discriminants $|d|\le x$ 
with $L(1,\chi_d) \le \zeta(2)/(e^{\gamma}\log_2 |d|)$.
\endproclaim 

\proclaim{Conjecture 2 (Montgomery-Vaughan)}  
 The proportion of fundamental discriminants 
$|d|\le x$ with $L(1,\chi_d) \ge e^{\gamma} (\log_2 |d| 
+\log_3 |d|)$ is $>x^{\theta}$ and $<x^{\Theta}$ where 
$0 <\theta <\Theta <1$.  Similar estimates apply to the frequency 
with which $L(1,\chi_d) \le \zeta(2)e^{-\gamma} (\log_2 |d| 
+\log_3 |d|)^{-1}$.  
\endproclaim

\proclaim{Conjecture 3 (Montgomery-Vaughan)}  
For any $\epsilon >0$ there are 
only finitely many $d$ with $L(1,\chi_d) > e^{\gamma} (\log_2 |d| 
+(1+\epsilon)\log_3 |d|)$, or with 
$L(1,\chi_d) \le \zeta(2)e^{-\gamma} (\log_2 |d| 
+(1+\epsilon)\log_3 |d|)^{-1}$.
\endproclaim

Notice that Conjecture 3 implies that 
the true nature of extreme values of $L(1,\chi_d)$ is 
given by Chowla's omega-results rather than the GRH bounds 
of Littlewood.  

The idea of comparing the distribution of values of 
$L(1,\chi_d)$ to a random model is quite old.  Chowla and P. Erd{\H o}s 
[3] showed first that for a fixed $\tau$ the proportion of 
$|d|\le x$ with $L(1,\chi_d) \ge e^{\gamma} \tau$ tends 
to a limit as $x \to \infty$.  Further this distribution function 
was shown to be continuous.  Elliott [6] pursued this further and 
showed that $L(1,\chi_d)$ possesses a smooth distribution function, and 
he also obtained expressions for its characteristic function 
(Fourier transform).  Indeed Elliott established such 
results for the distribution of $L(s,\chi_d)$ 
at any point $s$ with Re$(s)>1/2$, and he also allows for $d$ to be 
restricted to prime discriminants. 
Elliott's results also allow for $\tau$ to grow slowly in terms of $x$; 
essentially his methods 
show that the distribution of values of $L(1,\chi_d)$ approximate 
the distribution of random Euler products for $\tau$ throughout the range 
$(1/\log_3 x, \log_3 x)$.  But these results are not sufficiently 
uniform to approach the above Conjectures of Montgomery and Vaughan.

We now describe our results, starting with the probablistic 
model we shall use to study $L(1,\chi_d)$.  
For primes $p$ let $X(p)$ denote independent random variables taking the 
values $1$ with probability $p/(2(p+1))$,
$0$ with probability $1/(p+1)$, and $-1$ with probability $p/(2(p+1))$.
The reason for this choice of probabilities over the simpler 
$\pm 1$ with probability $1/2$ is 
as follows:  For odd primes $p$, fundamental discriminants $d$ 
are constrained to lie in one of $p^2-1$ residue classes $\pmod {p^2}$, 
and for $p-1$ of these the character value is $0$, and the remaining 
$p(p-1)$ residue classes split equally into $\pm 1$ values.  For 
the prime $2$ note that fundamental discriminants 
lie in the residue classes $1,5,8,9,12,13 \pmod {16}$ and the 
values $0$, $\pm 1$ occur equally often.  
We extend $X$ multiplicatively to all integers $n$: that is 
set $X(n)=\prod_{p^{\alpha} \parallel n} X(p)^{\alpha}$. 
We wish to compare the distribution of values of $L(1,\chi_d)$ with 
the disrtibution of values of the random Euler products $L(1,X):= 
\prod_p (1-X(p)/p)^{-1}$ (these products converge with probability $1$).

Before describing how well this model approximates the 
distribution of $L(1,\chi_d)$, it is helpful to gain an 
understanding of the distribution of $L(1,X)$.  To this end, 
we define
$$
\Phi(\tau) = \text{Prob}(L(1,X) \ge e^{\gamma} \tau), 
$$
and 
$$
\Psi(\tau) = \text{Prob} (L(1,X) \le \tfrac{\pi^2}{6} 
\tfrac{1}{e^{\gamma }\tau} ). 
$$
By studying the characteristic function of $\log L(1,X)$ 
(which may be shown to decay rapidly) one can see that 
$\Phi(\tau)$ and $\Psi(\tau)$ are 
smooth functions.  This is implicit in Elliott [6], and 
follows also from our subsequent work.  
From the work of Montgomery and Vaughan [15] we 
obtain that $\Phi(\tau)$ and $\Psi(\tau)$ 
decay double exponentially as $\tau\to \infty$:  precisely, 
there exist constants $C$ and $c$ such that 
$$
\exp\Big( - C\frac{e^{\tau}}{\tau}\Big) \le \Phi(\tau) 
\le \exp\Big(-c\frac{e^{\tau}}{\tau}\Big),
$$ 
and similarly for $\Psi(\tau)$.  In Section 3
we shall analyse these functions closely, expressing them in 
terms of the moments of $L(1,X)$ (see Theorem 3.1 below).  
From our work we obtain the following useful estimates 
for $\Phi(\tau)$ and $\Psi(\tau)$, which improve upon 
Montgomery and Vaughan's estimates. 

\proclaim{Proposition 1}   For large $\tau$ we have 
$$
\Phi(\tau) = \exp\Big( -\frac{e^{\tau -C_1}}{\tau} + O\Big(
\frac{e^{\tau}}{\tau^2}\Big)\Big), 
$$
where 
$$
C_1 := \int_{0}^1 \tanh y \frac{dy}{y} + \int_1^\infty (\tanh y-1) 
\frac{dy}{y} = 0.8187\ldots.  
$$ 
The same asymptotic holds also for $\Psi(\tau)$. Further 
if $0\le \lam \le e^{-\tau}$ then
$$
\Phi(\tau e^{-\lam}) 
= \Phi(\tau)(1+ O(\lam e^{\tau})), 
\qquad \text{and } \qquad
\Psi(\tau e^{-\lam}) 
= \Psi(\tau)(1+ O(\lam e^{\tau})).  
$$
\endproclaim

Below we shall let $\sumf$ indicate that the sum is over 
fundamental discriminants.  
A standard argument shows that there are $ \frac{6}{\pi^2} x
+O(x^{\frac 12+\epsilon})$ 
fundamental discriminants $d$ with $|d|\le x$. 
Define ${\Phi_x}(\tau)$ to be the proportion of the fundamental
discriminants $d$ with $|d|\leq x$, for which 
$L(1,\chi_d)>e^{\gamma}\tau$; that is
$$
{\Phi_x}(\tau):= \Big( \sumf\Sb|d|\le x\\ L(1,\chi_d) > e^{\gamma} 
\tau\endSb 1 \Big)\Big/
\Big( \sumf\Sb|d|\le x\endSb 1 \Big)
$$
(and similarly define $\Psi_x(\tau)$).
We would like to 
compare this with $\Phi (\tau)$ (and analogously 
to compare the frequency of small values with $\Psi$).  Notice 
that the viable range for such a correspondence is 
$\tau \le \log_2 x + \log_3 x + C_1 + o(1)$.  
Proposition 1 shows that at this juncture the probabilities 
$\Phi$ and $\Psi$ become smaller than $1/x$.

\proclaim{Theorem 1}  Let $x$ be large. 
Uniformly in the region $\tau \le R(x)$ 
we have 
$$
{\Phi_x}(\tau)= {\Phi}(\tau) \Big( 1 + O\Big( \frac{1}{(\log x)^5} + 
e^{\tau -R(x)}\Big)\Big), 
$$
and uniformly in $\tau \le R(x)+ \log_3 x$ we have 
$$
{\Phi_x}(\tau) = {\Phi}(\tau) (\log x)^{O(1)}.
$$
Here we may choose $R(x)= \log_2 x -2\log_3 x+\log_4 x -20$ unconditionally, 
and $R(x) = \log_2 x -\log_3 x  -20$ if the GRH is assumed.  
Analogous results hold replacing $\Phi$ with 
$\Psi$.
\endproclaim 

Our proof of Theorem 1 relies upon computing the mean moment of $L(1,\chi_d)^z$,
as we average over fundamental discriminants $d$ with $|d|\le x$, 
for complex numbers $z$ in a wide range.  
We will establish that these moments are very nearly equal to the 
corresponding moments of the random $L(1,X)$; that is to say the 
expectation ${\Bbb E}(L(1,X)^z)$.  Throughout the paper ${\Bbb E}(\cdot)$ 
stands for the expectation of the random variable in brackets. 
Suppose there were a character $\chi_d$ with $|d|\asymp x$ for which $L(s,\chi_d)$ 
has a bad Landau-Siegel zero.  In this case we could have $L(1,\chi_d)$ 
as small as $x^{-\epsilon}$ so that when $z$ is a negative real number 
the mean moment of $L(1,\chi_d)^z$ 
would be heavily affected by this particular character. Thus, short of 
proving the non-existence of Landau-Siegel zeros, we cannot hope for 
asymptotics for the moments as stated, except in a narrow range of values 
for $z$.  To circumvent this difficulty, we calculate instead moments of 
$L(1,\chi_d)$ after first omitting a sparse set of discriminants having 
Landau-Siegel zeros.  Precisely, we define for an appropriate
constant $c>0$ 
$$
{\Cal L} = \{ d : \qquad L(\beta,\chi_d) =0 \qquad \text{for some   }
 1-c/\log (e|d|) \le \beta <1\}. 
$$
We will refer to elements of ${\Cal L}$ as Landau-Siegel characters, 
discriminants, or moduli.  If $c$ is chosen appropriately, then 
it is known that there is at most one element $d$ of ${\Cal L}$ between 
$x$ and $2x$ (see [4]) so that there are $\ll \log x$ elements of ${\Cal L}$ 
with $|d|\le x$.  

\proclaim{Theorem 2}  Uniformly in the region 
$|z| \le \log x/(500 (\log_2 x)^2)$ we have
$$
\sumf\Sb |d|\le x\\ d\notin {\Cal L} \endSb L(1,\chi_d)^z 
= \frac{6}{\pi^2} x {\Bbb E}(L(1,X)^z) + O\Big( x \exp\Big(
-\frac{\log x}{5\log_2 x}\Big)\Big).
$$
\endproclaim

We remark that when Re $z$ is positive then it is 
not necessary to omit elements of ${\Cal L}$ while calculating 
the moments of Theorem 2.  This follows because $L(1,\chi_d)$ is 
easily seen to be $\ll \log (e|d|)$ and so for positive Re $z$ the error 
in adding back elements of ${\Cal L}$ is $\ll \log x (\log x)^{\text{Re }z}$ 
which may be subsumed in the error term of the Theorem.  Previously 
A.F. Lavrik [10] had computed $\sumf_{|d|\le x} L(1,\chi_d)^{2k}$ for integers 
$k \ll \sqrt{\log x}$, and we see that Theorem 2 goes substantially 
further than his result. Further if Re $z$ is negative, but $|\text{Re }z|$ 
is bounded, then again we need not exclude elements of ${\Cal L}$.  This 
follows readily from Siegel's famous 
bound  $L(1,\chi_d) \gg_{\epsilon} |d|^{-\epsilon}$.

We may extend the range of applicability (in $|z|$) of Theorem 2 by 
excluding a larger (but still very thin) set of characters.  We 
describe this result next, which will be the main ingredient 
used to prove Theorem 1.  

\proclaim{Theorem 3}  Let ${\Cal E}$ denote a set of 
$\le \sqrt{x}$ exceptional discriminants $|d|\le x$.  
Uniformly in the region $|z|\le \log x \log_3 x/(e^{12} \log_2 x)$  
we have 
$$
\sumf\Sb |d|\le x\\ d\notin {\Cal E} \endSb L(1,\chi_d)^z 
= \frac{6}{\pi^2} x {\Bbb E}(L(1,X)^z) 
+ O\Big( x \frac{{\Bbb E}(L(1,X)^{\text{Re }z})}{(\log x)^9}\Big).
$$
If the GRH is true then the above asymptotic 
holds uniformly in the larger region $|z| \le 10^{-3}\log x$.  
\endproclaim

One can show (by modifying Lemma 3.2 below) that the main term in 
Theorem 3 dominates the error term if and only if $|\text{Im } z| \ll 
\{(1+|\text{Re }z|)\log (2+|\text{Re }z|)\}^{\frac 12} \log_2 x$.  
In contrast the main term in Theorem 2 dominates the error term 
throughout the region $|z|\le \log x/(500 (\log_2 x)^2)$.  

Even assuming GRH the range of validity of Theorem 1 is insufficient
to penetrate the Conjectures of Montgomery and Vaughan.  However 
note that these conjectures ask only that the frequencies of large 
and small values decay ``double exponentially'' and not for 
the more precise information provided by Theorem 1.  We now set ourselves 
the intermediate problem of determining when ${\Phi_x}(\tau), {\Psi_x}(\tau) = 
\exp(-(1+o(1))e^{\tau-C_1}/\tau)$.  Exploiting a wonderful result of S.W. Graham and 
C.J. Ringrose [8] on character sums to smooth moduli (see Lemma 4.2 below), 
we settle this problem in a range wider than required for Conjecture 1.

\proclaim{Theorem 4}  Let $x$ be large 
and let $\log_2 x \ge A \ge e$ be a real number.  Uniformly in the 
range $\tau \le R_1(x) -\log_2 A$ we have 
$$
{\Phi_x}(\tau) = \exp\Big(-\frac{e^{\tau-C_1}}{\tau}\Big(1+O\Big(\frac 1A +
\frac 1{\tau}\Big)\Big)\Big), 
$$
and the same asymptotic holds ${\Psi_x}(\tau)$.  
Here we may take $R_1(x)= \log_2 x+\log_4 x -20$ unconditionally, and 
$R_1(x) = \log_2 x +\log_3 x -20$ if the GRH is true.  
\endproclaim  

Theorem 4 clearly implies a stronger version of Montgomery and Vaughan's 
Conjecture 1.  Notice also that assuming GRH it implies part of Conjecture 2:
namely that the number of $|d| \le x$ with $L(1,\chi_d) \ge e^{\gamma}
(\log_2 |d| +\log_3 |d|)$ is $\le x^{\Theta}$ for some $\Theta <1$, and 
similarly for small values.  In [13]  Montgomery established 
that if the GRH is true then there are infinitely many primes $p$ such 
that the least quadratic non-residue $\pmod p$ is $\gg \log p \log \log p$.  
We adapt his idea to examine extreme values of $L(1,\chi)$ under GRH, which 
strengthens Theorem 4, but just fails to prove the other half of 
Conjecture 2.

\proclaim{Theorem 5a} Assume GRH. For any $\epsilon>0$, and all
large $x$, there are $\gg x^{\frac 12}$ primes $q\le x$ such that
$$
L(1,(\tfrac{\cdot}{q})) \ge e^{\gamma}(\log_2 q +\log_3 q 
- \log(2\log 2) -\epsilon), 
$$
and $\gg x^{\frac 12}$ primes $q$ such that 
$$
{L(1,(\tfrac{\cdot}{q}))} \le 
\frac{\zeta(2)}{e^{\gamma}} (\log_2 q +\log_3 q
- \log(2\log 2) -\epsilon)^{-1}. 
$$
\endproclaim

Note that $\log (2\log 2) = 0.3266\ldots$ so that Theorem 5a 
comes very close to exhibiting the extreme values required in 
Conjecture 2.  We may ask for the extreme values of $L(1,\chi_d)$ that 
may be obtained unconditionally: that is, for refinements of 
Chowla's results which would approach the extreme values 
predicted by the Conjectures above.  As far as small values are concerned 
we are unable to go further than the extreme values guaranteed 
by Theorem 1.  However by a judicious use of the pigeonhole principle 
we are able to exhibit large values of $L(1,\chi_d)$ which are 
nearly as good as the conjectured truth.

\proclaim{Theorem 5b}  For large $x$ there are at least $x^{\frac 1{10}}$ 
square-free integers $d \le x$ such that 
$$
L(1,(\tfrac{\cdot}{d})) \ge e^{\gamma}(\log_2 x+ \log_3 x -\log_4 x -10).
$$
\endproclaim

In summary our results find excellent agreement between the distribution of 
values of $L(1,\chi_d)$ and the predictions of the probabilistic 
model.  We find (especially on GRH) that the predictions 
hold true in almost the entire viable range, and this leads us to 
believe that the extreme values of $L(1,\chi_d)$ behave like 
Chowla's omega results.  

If the asymptotic formula of 
Theorem 1 holds to the edge of the viable range then perhaps
$$
\max\Sb d\ \text{\rm fundamental}\\ |d|\leq x \endSb  \ L(1,\chi_d) =
e^\gamma ( \log_2 x+ \log_3 x +C_1 +o(1) ) .
$$
It is also plausible that the distribution function changes nature just beyond
the range given in Theorem 1 (which is why our methods 
do not give good results there),
and that the maximum value is represented by a 
slightly different function. Nonetheless
given the extraordinary decay of 
$\Phi(\tau)$ that we have detected in this range
we conjecture that, at worst, the above estimate is 
true with a slightly different constant.

\head 2.  Preparatory Lemmas \endhead 

\noindent In this section we collect together some preliminary 
results which will be useful in our subsequent work.  The ideas in 
this section are standard and classical.

We will show in Proposition 2.2 below that with few exceptions $L(1,\chi)$ 
may be approximated by the short Euler product 
$L(1,\chi;y)$ where, throughout the paper, we let $L(1,\chi;y):= \prod_{p\le y} 
(1-\chi(p)/p)^{-1}$. 
The primary ingredient is the following classical lemma,
which is proved in Lemmas 8.1 and 8.2 of [9].

\proclaim{Lemma 2.1}  Let $s=\sigma+it$ with $|t|\le 3q$, and let $y\ge 2$ 
be a real number.  Let $\frac 12\le \sigma_0 < \sigma$, and 
suppose that the rectangle $\{ z: \ \sigma_0 < \text{Re}(z) \le 1, \ \ 
|\text{Im}(z)-t|\le y+ 3\}$ contains no zeros of $L(z,\chi)$.  
Then 
$$
|\log L(s,\chi)| \ll \frac{\log q}{\sigma-\sigma_0}.
$$
Further, putting $\sigma_1= \min(\sigma_0 +\frac{1}{\log y}, 
\frac{\sigma+\sigma_0}{2})$, 
$$
\log L(s,\chi)  = \sum_{n=2}^{y} \frac{\Lambda(n)\chi(n)}{n^s \log n} + 
O\biggl(\frac{\log q}{(\sigma_1-\sigma_0)^2} 
y^{{\sigma_1-\sigma}}\biggr).
$$
\endproclaim

If the GRH holds then we may apply Lemma 2.1 with $\sigma_0=\frac 12$. 
We now take $y=(\log q)^2(\log_2 q)^6$ to obtain 
$$
\log L(1,\chi) =\sum_{n=2}^{y} \frac{\Lambda(n)\chi(n)}{n\log n} +
O\left( \frac{1}{ \log\log q}\right).
$$ 
Using the prime number theorem to estimate the contribution
of the primes between $\log^2q$ and $y$, we deduce that
$$
L(1,\chi) = \prod_{p\leq \log^2q} \left( 1- \frac{\chi(p)}{p}
\right)^{-1} \left\{ 1+ 
O\left( \frac{\log_3 q}{ \log_2 q}\right) \right\},
$$
which gives Littlewood's result (1.1).  

Using this lemma together with the large sieve and 
zero density results we obtain the following result,  
essentially due to Elliott [5].

\proclaim{Proposition 2.2}  Let $Q$ be large, and 
let $\log_2 Q \ge A\ge 1$ be a real number.  Then 
for all but at most $Q^{2/A}$ primitive characters $\chi \pmod q$ with 
$q\le Q$ we have for $y\le Q/2$ 
$$
L(1,\chi) = \prod_{p\le y} \Big(1-\frac{\chi(p)}{p}\Big)^{-1} 
\Big(1+O\Big(\frac{A^2 \log Q }{y^{\frac 1{4A}}} \Big)\Big). \tag{2.1}
$$
Further  
$$
L(1,\chi) = \prod_{p\le (\log Q)^A} \Big(1-\frac{\chi(p)}{p}\Big)^{-1} 
\Big(1+ O\Big(\frac{1}{\log \log Q}\Big)\Big)
$$
holds for all but at most $Q^{2/A + 5\log_3 Q/\log_2 Q}$ primitive 
characters $\chi \pmod q$ with $q\le Q$.  
\endproclaim

\demo{Proof}  From a standard zero density result 
(see Theorem 20 of E. Bombieri [1]) we know that there 
are fewer than $Q^{6(1-\alpha)} (\log Q)^B$ primitive 
characters with conductor below $Q$ having a 
zero in the rectange $1\ge \text{Re }(s)\ge \alpha$, $|\text{Im }(s)|\le Q$.  
Here $B$ is some absolute constant.  Thus appealing  to Lemma 2.1 
with  $s=1$,  and $\sigma_0 =1 -\frac{1}{4A}$ 
we obtain (2.1) for all but at most $Q^{2/A}$ primitive characters 
$\chi \pmod q$ with $q\le Q$.  

We now show that 
$$\sum_{(\log Q)^{A} \le p \le (\log Q)^{8A}} \frac{\chi(p)}p 
=O\Big(\frac{1}{\log_2 Q}\Big)
$$
for all but $Q^{2/A +5\log_3 Q/\log_2 Q}$ characters 
with conductor below $Q$, which when combined with (2.1) with 
$y=(\log Q)^{8A}$  gives the Proposition.  
To prove this we shall use the large sieve in the following 
form:
$$
\sum_{q\le Q} \qquad 
\sumstar_{\chi \pmod q} \Big| \sum_{m\le M} a(m) \chi(m)\Big|^2 
\ll (Q^2+M ) \sum_{m\le M} |a(m)|^2, \tag{2.2} 
$$ 
where the $\sumstar$ is over primitive characters $\chi$, and the $a(m)$ 
are arbitrary complex numbers.  

For $0\le j \le J:=[7A \log_2 Q/\log 2]$  put $z_j= 2^j (\log Q)^A$ and 
put $z_{J+1} = (\log Q)^{8A}$. 
Choose $k=[2\log Q/(A\log_2 Q)]+1$ so that $M_j:=z_{j+1}^{k} \ge Q^2$.  
Define 
$$
\sum_{M_j/2^{k} \le m\le M_j} a_j(m) \frac{\chi(m)}{m} = 
\Big(\sum_{z_j <p \le z_{j+1}} \frac{\chi(p)}p\Big)^{k}.
$$
Note that $0\le a_j(m) \le k!$ 
and $\sum_{m} a_j(m) = (\pi(z_{j+1})-\pi(z_j))^{k} 
\le (2z_{j+1}/(A\log_2 Q))^{k}$.  
Appealing to (2.2) we see that 
$$
\align
\sum_{q\le Q}\qquad \sumstar_{\chi \pmod q} \Big|
\sum_{z_j <p \le z_{j+1}} \frac{\chi(p)}p\Big|^{2k}
&\ll z_{j+1}^{k} \sum_{M_j/2^{k} \le m\le M_j} \frac{|a_j(m)|^2}{m^2} 
\ll \frac{k! 2^{2k} }{z_{j+1}^{k}} \Big(\frac{2z_{j+1}}{A\log_2 Q}\Big)^{k} 
\\
&\ll \Big( \frac{4k}{A \log_2 Q}\Big)^{k},
\\
\endalign
$$
using Stirling's formula.  We deduce that 
$$
\Big|\sum_{z_j\le p\le z_{j+1} } \frac{\chi(p)}{p}\Big| \le 
\frac{1}{A (\log_2 Q)^2}
$$
for all but at most $(10 \log Q (\log_2 Q)^2)^k$ primitive characters with 
conductor below $Q$, and the Proposition follows at once.

\enddemo

Using Lemma 2.1, we may also derive the following approximation 
to $L(1,\chi)^z$, provided $\chi$ is not a Landau-Siegel character. 

\proclaim{Lemma 2.3}  Let $\chi$ be a non-principal character 
$\pmod q$ which is either complex, or real and primitive with its  
discriminant not in ${\Cal L}$.  
Let $z$ be any complex number with $|z| \le (\log q)^2$ and
let $Z \ge \exp((\log q)^{10})$ be a real number.  
Then 
$$
L(1,\chi)^z = \sum_{n=1}^{\infty} \chi(n) \frac{d_z(n)}{n} e^{-n/Z} 
+ O\left( \frac 1q \right).
$$
\endproclaim

\demo{Proof}  Since $\frac{1}{2\pi i} \int_{1-i\infty}^{1+i\infty} y^s 
\Gamma(s) ds = e^{-1/y}$, we have
$$
\frac{1}{2\pi i} \int_{1-i\infty}^{1+i\infty} L(1+s,\chi)^z Z^s \Gamma(s) ds 
= \sum_{n=1}^{\infty} \frac{d_z(n)}{n} \chi(n) e^{-n/Z}. \tag{2.3}
$$
We shift the line of integration to the contour $s=-C(t)+it$ where
$C(t): = -c/(2\log (q(|t|+2)))$, for an appropriately small constant $c>0$.  
From our assumption on $\chi$ we may choose $c$ such that we 
encounter no zeros of $L(s,\chi)$ while shifting contours.  
We encounter a pole at $s=0$ which leaves the
residue $L(1,\chi)^z$.
Applying Lemma 2.1 with $\sigma_0=1-4C(t)/3$ and $y=2$ gives 
$|\log L(s,\chi)|\ll 1/C(t)^2$ and so 
 the left side of (2.3) equals $L(1,\chi)^z + R$ where 
$$
R \ll \int_{-\infty}^{\infty} Z^{-C(t)} 
e^{O(|z|/C(t)^2)} | \Gamma(-C(t) +it)| dt \ll \frac{1}{q},
$$
by Stirling's formula.
\enddemo

We end this section by collecting several inequalities 
for values of $d_z(n)$: the $z$-th divisor function. 
Recall that $d_z(n)$ is a multiplicative function 
given on prime powers by $d_z(p^a)= \Gamma(z+a)/(\Gamma(z)a!)$, and 
that $\sum_{n=1}^{\infty} d_z(n)/n^s =\zeta(s)^z$ for Re $s>1$. 
 
Note that $|d_z(n)|\le d_{|z|}(n)$.  
For real numbers $k\ge 1$ we observe that $d_k(mn)\le d_k(m)d_k(n)$.  
For any positive integers $a,b,n$ we have $d_a(n)d_b(n)\leq d_{a+b}(n)$.
Further for any complex number $z$ and a real number $\beta$ we 
have $|d_z(n)|^{\beta} \le d_{|z|^{\beta}}(n)$.  We also record that 
$|d(n)d_z(n^2)| \le d_{2|z|+2}(n)^2$ and that $|d_z(n^2)| \le 
d_{(|z|+1)^2}(n)$.  All of these inequalities may be shown by 
by first proving them for prime powers (by induction on the exponent), 
and then using multiplicitivity to deduce them for all $n\geq 1$. 

Lastly we note that for positive integers $k$ and real numbers $x\ge 2$ we 
have 
$$
\sum_{n\le x} \frac{d_k(n)}{n} \le \Big(\sum_{n\le x} \frac 1n\Big)^k 
\le (\log 3x)^k, 
$$
and since $d_k(n) e^{-n/x} \le e^{k/x} \sum_{a_1 \ldots a_k =n}
e^{-(a_1+\ldots +a_k)/x}$ that 
$$
\sum_{n=1}^{\infty} \frac{d_k(n)}{n} e^{-n/x} \le \left( e^{1/x} 
\sum_{a=1}^{\infty} \frac{e^{-a/x}}{a} \right)^k \le (\log 3x)^k.
\tag{2.4}
$$

\head 3. Random Euler Products and their distribution \endhead

\noindent We put $L(1,X;y)= \prod_{p\le y} (1-X(p)/p)^{-1}$. 
Since ${\Bbb E}((\sum_{p> y} X(p)/p)^2)\asymp  \sum_{p>y}\frac{1}{p^2} 
\asymp 1/(y\log y)$, we 
see that with probability $1$, $L(1,X;y)$ converges to $L(1,X)$ as 
$y\to \infty$.  In this section we investigate the distribution of 
the random Euler products $L(1,X;y)$.  Letting $y=\infty$ in our 
results gives information on the distribution of $L(1,X)$.   
Analogously to the definitions of $\Phi(\tau)$ and $\Psi(\tau)$ we define 
$$
{\Phi}(\tau;y) = \text{Prob}(L(1,X;y)\ge e^{\gamma} \tau), 
\ \ \  \text{and    }  {\Psi}(\tau;y)= \text{Prob}(L(1,X;y)\le 
\tfrac {\pi^2}{6}\tfrac 1{e^{\gamma} \tau}). 
$$

To facilitate our discussion we define for positive real numbers $k$ 
$$
I(k;y):=  -\sum_{p\le y} \log \Big(1-\frac{1}{p}\Big) \tanh \L{k}{p},
\ \ \ \text{and} \ \ \ I(k):= - \sum_{p}  
\log \Big(1-\frac{1}{p}\Big) \tanh \L{k}{p}. 
$$
For fixed $y$ plainly $I(k;y)$ is an increasing function of $k$, and for 
fixed $k$ it is an increasing function of $y$.  
Further writing $g(t)=\tanh(t)$ if $t\le 1$, and $g(t)=\tanh(t)-1$ if 
$t>1$ we have that 
$$
I(k;y) = -\sum_{p\le \min(k,y)} \log\Big(1-\frac 1p \Big) 
- \sum_{p\le y} \log \Big(1-\frac{1}{p}\Big) g \L{k}{p},
$$
and now using 
the prime number theorem and partial summation we 
see that for large $k$ and $y$, 
$$
I(k;y) 
= \log_2 \min(k,y) +\gamma + \frac{C_1(k/y)}{\log k} 
+ O\Big( \frac{1}{\log^2 k} \Big), 
\tag{3.1}
$$
where 
$$
C_1(x) := \int_{x}^{\infty} \frac{g(t)}{t} dt. 
$$
Note that $C_1(0)= C_1=0.8187\ldots$.  
Define next
$$
R(k;y) =\sum_{p\le y} \frac{1}{p^2 \cosh^2 (k/p)}, \ \ \ 
\text{and    }  R(k)= \sum_p \frac{1}{p^2 \cosh^2(k/p)}. 
$$
Again appealing to the prime number theorem and 
partial summation we see that for large $k$ and $y$
$$
R(k;y) = \frac{1}{k\log k}\int_{k/y}^{\infty} \frac{dx}{\cosh^2 x} 
+ O\Big(\frac{1}{k\log^2 k}\Big) 
= \frac{1-\tanh(k/y)}{k\log k} + O\Big(\frac{1}{k\log^2 k}\Big).  \tag{3.2}
$$
We are now in a position to state 
our main result, which obtains an asymptotic 
formula for $\Phi(\tau;y)$ (and $\Psi(\tau;y)$) 
in terms of appropriate moments ${\Bbb E}(L(1,X;y)^k)$.  Naturally, 
letting $y\to \infty$ below furnishes asymptotics for $\Phi(\tau)$
or $\Psi(\tau)$.

\proclaim{Theorem 3.1}  Let $y$ be a large real number, and 
let $\tau$ be large and below $\log y- 1$.  Let 
$k=k_{\tau,y}$ denote the unique real number such that $I(k;y) = 
\gamma +\log \tau$.  
Then 
$$
\Phi(\tau;y) = \frac{{\Bbb E}(L(1,X;y)^{k})}
{k (e^{\gamma}\tau)^{k}} \frac{1}{\sqrt{2\pi R(k;y)}}
\Big(1+ O\Big(\frac{\log^2 k}{k^{\frac 14}}\Big)\Big), 
$$
and 
$$
\Psi(\tau;y) = \frac{{\Bbb E}(L(1,X;y)^{-k})}
{k(\frac 6{\pi^2}e^{\gamma }\tau)^k}  \frac{1}{\sqrt{2\pi R(k;y)}}
\Big(1+ O\Big(\frac{\log^2 k}{k^{\frac 14}}\Big)\Big).
$$
Further if $0 \le \lam \le e^{-\tau}$ then 
$$
{\Phi}(\tau e^{-\lam};y) -{\Phi}(\tau;y) 
\ll {\Phi}(\tau;y) \Big( \lam e^{\tau} + \frac{e^{\frac 34 \tau} 
\log y}{y} \Big), 
$$
and 
$$
{\Psi}(\tau e^{-\lam};y) -{\Psi}(\tau;y) 
\ll {\Psi}(\tau;y) \Big( \lam e^{\tau} + \frac{e^{\frac 34 \tau} 
\log y}{y} \Big).
$$
\endproclaim

We shall focus only on proving the results for $\Phi(\tau;y)$, 
the argument for ${\Psi}$ requires only some minor 
adjustments.  
For the rest of this section we write $s=k+it$, and we denote 
${\Bbb E}((1-X(p)/p)^{-s})$ by $E_p(s)$.  Clearly $|E_p(s)|\le E_p(k)$ 
always.

\proclaim{Lemma 3.2}  Let $s=k+it$, and let $k$ be large.  
If $p> k/4  $ then we have for some positive constant $c_0$ 
$$
|E_p(s)| \le E_p(k)  \exp\Big(-c_0\Big(1-\cos\Big(t\log\Big(\frac{p+1}{p-1}
\Big)\Big)\Big)\Big). 
$$
Further if $k \le 2y$ then there exists a positive constant $c$ such 
that  
$$
\Big|\frac{{\Bbb E}(L(1,X;y)^s)}{{\Bbb E}(L(1,X;y)^k)}\Big| 
\le 
\cases \exp(-c\frac{t^2}{k\log k}) &\text{if  } |t|\le k/4\\
\exp(-c\frac{|t|}{\log |t|}) &\text{if  } k/4 \le |t|\le y/2\\
\exp(-c\frac{y}{\log^3  y}) &\text{if  } y/2\le |t|\le y^{2}/\log^2 y.\\
\endcases
$$
\endproclaim

Note that letting $y\to \infty$ we get that 
${\Bbb E}(L(1,X)^s)/{\Bbb E}(L(1,X)^k) \ll \exp(-c \frac{t^2}{k\log k})$ 
if $|t|\le k$ and $\ll \exp(-c \frac{|t|}{\log |t|})$ if $|t| >k$.  

\demo{Proof}  Note that for $r_1$, $r_2$ and $r_3$ positive we 
have $|r_1+r_2e^{i\theta_2} + r_3 e^{i\theta_3}|^2 \le (r_1+r_2+r_3)^2 -
2r_1r_2 (1-\cos \theta_3)$, so that 
$|r_1+r_2e^{i\theta_2} 
+r_3 e^{i\theta_3}| 
\le (r_1+r_2+r_3) \exp(-\frac{r_1r_3(1-\cos \theta_3)}{(r_1+r_2+r_3)^2})$.  
We apply this with $r_1 =\frac{p}{2(p+1)} (1-1/p)^{-k}$, $r_2=1/(p+1)$,  
$r_3 = \frac{p}{2(p+1)} (1+1/p)^{-k}$, and $\theta_2= t\log (1-1/p)$ and 
$\theta_3= t \log \L{p-1}{p+1}$.  Since $p>k/4$ the first statement of 
the lemma follows.

To prove the second part, note that our desired ratio is 
$$
\le \prod_{y\ge p \ge k/4} \exp\Big(-c_0 \Big(1-\cos\Big(t
\log\Big(\frac{p+1}{p-1}\Big) \Big)\Big)\Big). 
$$
If $|t|\le k/4$ then for the primes $p$ in $(k/4,k/2)$ we have 
that $|t \log ((p+1)/(p-1))| \sim 2|t|/p$ lies between $\sim 4|t|/k$ 
and $\sim 8|t|/k$ so that $1-\cos(t \log((p+1)/(p-1))) \gg |t|^2/k^2$.  
This gives the bound of the Lemma in this case.  Next if 
$k/4 \le |t|\le y/2$ then we apply the above argument 
with the primes in $(|t|,2|t|)$ getting the desired bound.  

Lastly suppose that $y/2 \le |t|\le y^2/\log y$.  Here we let 
$\delta := 10^{-6}/\log y$, and divide the interval $(y/2,y)$ into 
intervals of length $\delta y^2/|t|$ (with the last interval possibly 
being shorter).  There are $\sim |t|/(2y\delta)$ such intervals.  Call 
an interval good if $\cos(t\log(\frac{p+1}{p-1})) \le \cos(\delta/10)$ 
for all primes $p$ in that interval, and bad if otherwise.   
There are at most $3|t|/y$ bad intervals, and by the Brun-Titchmarsh theorem 
each bad interval contains at most $3\delta y^2/(|t| \log (\delta y^2/|t|))
\le y^2/(200 |t| \log y)$ primes.  Thus there are at least $y/3\log y$ 
primes in good intervals, and the Lemma follows in this final case.

\enddemo

In addition to Lemma 3.2 we need a result comparing 
${\Bbb E}(L(1,X;y)^{s})$ with ${\Bbb E}(L(1,X;y)^k)$ for 
relatively small values of $t$. First observe that 
$$
{\Bbb E} \Big(\Big(1-\frac{X(p)}p\Big)^{-k}\Big) = {\Bbb E}\Big( 
\exp\Big(k\frac{X(p)}{p}\Big)\Big) \exp\Big(O\Big(\frac{k}{p^2}\Big)\Big)
=
\cosh \L{k}{p} \exp\Big(O\Big(\frac k{p^2}\Big)\Big), \tag{3.3} 
$$ 
where the last estimate follows since ${\Bbb E}(\exp(kX(p)/p)) 
= (p\cosh (k/p)+1)/(p+1)$ which equals $\cosh(k/p)(1+O(1/p))$ 
if $p\le k$ and equals $\cosh(k/p)+O(k^2/p^3) = \cosh(k/p)(1+O(k^2/p^3))$ 
if $p>k$.  Further note that if $p\le k/\log k$ then 
$$
{\Bbb E}\Big(\Big(1-\frac{X(p)}p\Big)^{-k}\Big) 
= \Big(\frac{p}{2(p+1)}\Big) \Big(1-\frac{1}{p}\Big)^{- k} 
\Big( 1+ O\Big(\frac{e^{-k/p }}{p}\Big)\Big).
\tag{3.4}
$$ 
Note that,  for all primes $p$ we have
$$
\align
E_p(s) &= 
\Big(1-\frac 1p\Big)^{-it} {\Bbb E}\Big( \Big(1-\frac {X(p)}{p}\Big)^{-k} 
\Big(\frac{p-1}{p-X(p)}\Big)^{it}\Big) .\\
\endalign
$$
Now
$$
\left(\frac{p-1}{p-X(p)}\right)^{it} = 1-it\frac{(1-X(p))}p - t^2 \frac{(1-X(p))^2}{2p^2} + 
O\left((1-X(p)) \left(   \frac{|t|}{p^2} + \frac{|t|^3}{p^3}\right)\right) :
$$
this is trivial if $X(p)=1$, otherwise it follows from the Taylor
expansion if $p>|t|$, and as $|(p-1)/(p-X(p))^{it}| = 1$ if $p\leq |t|$.
 By (3.4) we have that 
${\Bbb E}((1-X(p))(1-X(p)/p)^{-k}) = (1+1/p)^{-k} +O(1/p)$ and 
${\Bbb E}((1-X(p))^2 (1-X(p)/p)^{-k}) = 2(1+1/p)^{-k} +O(1/p)$,
 and so we deduce that 
$$
E_p(s) 
=\Big(1-\frac 1p\Big)^{-it} \Big( E_p(k) 
- \frac{it}{p} \left( 1+ \frac 1p\right)^{-k} 
-\frac{t^2}{p^2} \left( 1+ \frac 1p\right)^{-k} 
+O\Big(\frac{|t|}{p^2} + \frac{|t|^2}{p^3} + \frac{|t|^3}{p^3}\Big)\Big).
$$

If $p\leq k/3\log k$ then $E_p(k)\geq \frac 13 (1-1/p)^{-k}\gg k^3$
and so the above becomes, using (3.3),
$$
E_p(s) 
=\Big(1-\frac 1p\Big)^{-it} E_p(k)
 \left( 1 +O\Big(\frac{|t|+|t|^3}{k^3}\Big)\right).
$$
If $p\geq k/3\log k$ then $(1+1/p)^{-k}=e^{-k/p}+O(1/p)$. If 
$k/3\log k\leq p\leq |t|$ then 
$$
E_p(s) 
=\Big(1-\frac 1p\Big)^{-it} E_p(k)
 \left( 1 +O\Big(\frac{|t|^3}{p^3} e^{-k/p}\Big)\right).
$$
If $p\geq k/3\log k$ and $p\geq |t|$ then the above becomes, using (3.3),
$$ 
\log \frac{E_p(s)}{E_p(k)}=  -it \log \Big(1-\frac 1p\Big) \tanh 
\Big(\frac kp\Big) -\frac{t^2}{2} \frac{1}{p^2 \cosh^2 (k/p)} 
+ O\Big(\Big(\frac{|t|}{p^2}+  \frac{|t|^3}{p^3}\Big) 
e^{-k/p}\Big).
$$
It follows that, for $|t|\leq k^{2/3}$,
$$
{\Bbb E}(L(1,X;y)^s) = {\Bbb E}(L(1,X;y)^k) \exp\Big(it I(k;y) - \frac{t^2}{2} 
R(k;y)+ O\Big(\frac{|t|}{k\log k}+\frac{|t|^3}{k^2 \log k}\Big)\Big). 
\tag{3.5}
$$
We are now in a position to prove Theorem 3.1.  

\demo{Proof of Theorem 3.1}  Let $\lam >0$ be real.  
For $y>0$ and $c>0$ we have by Perron's formula 
$$
\align
\frac{1}{2\pi i} \int_{c-i\infty}^{c+i\infty} y^z \Big(\frac{e^{\lam z}-1}
{\lam z}\Big) \frac{dz}{z} &= \frac{1}{\lam } \int_0^{\lam} 
\left( \frac{1}{2\pi i} \int_{c-i\infty}^{c+i\infty} (ye^{u})^z \frac{dz}{z} \right) du \\
&= \ \cases 1 &\text{if  } y>1 \\
1+\frac{\log y}\lam \in [0,1] &\text{if  } e^{-\lam} \le y \le 1\\
0 &\text{if  } y< e^{-\lam}.\\
\endcases \\
\endalign
$$
Also by Perron's formula we may see that 
$$
\frac{1}{2\pi i } \int_{c-i\infty}^{c+i\infty} y^z \Big(\frac{e^{\lam z}-1}
{\lam z}\Big) \Big(\frac{e^{\lam z}-e^{-\lam z}}{z}\Big) dz
$$
is always non-negative, and that it equals $1$ if $e^{-\lam} \le y\le 1$.  
Applying these identities we obtain that 
$$ 
\Phi(\tau;y) \le 
\frac{1}{2\pi i} \int_{k-i\infty}^{k+i\infty} {\Bbb E}(L(1,X;y)^s) 
(e^{\gamma}\tau)^{-s} \Big(\frac{e^{\lam s}-1}{\lam s}\Big) \frac {ds}{s}
\le \Phi(\tau e^{-\lam};y), \tag{3.6}
$$ 
and 
$$
\Phi(\tau e^{-\lam};y)-\Phi(\tau;y)
\le \frac{1}{2\pi i} \int_{k-i\infty}^{k+i\infty} {\Bbb E}(L(1,X;y)^s) 
(e^{\gamma}\tau)^{-s}\Big(\frac{e^{\lam s}-1}{\lam s}\Big) \Big(
\frac{e^{\lam s}-e^{-\lam s}}{s}\Big) ds. \tag{3.7}
$$

Now suppose that $\tau \le \log y -1$ is large.  We choose 
$k= k_{\tau,y}$ such that $I(k;y) = \gamma+ \log \tau$.  Observe that 
$k\le y$ since $I(k;y)$ is increasing in $k$ and by (3.1) $I(y;y) = 
\log_2 y +\gamma +C_1(1)/\log y+O(1/\log^2 y) > \log \tau +\gamma$ as
$C_1(1) = -0.09\ldots$.   Further (3.1) gives 
us that $k \asymp e^{\tau}$.  We now suppose below that $\lam \le e^{-\tau} 
\ll 1/k$, so that $|e^{\lam s}| =e^{\lam k} \ll 1$.

Consider first the integral in (3.7).  We split the integral into the 
part when $|t|\le y^2/\log^2 y$, and when $|t|> y^2/\log^2 y$.  The
contribution of the second segment of the integral is 
$$
\ll \frac{{\Bbb E}(L(1,X;y)^{k})}{(e^{\gamma} \tau)^k} 
\int_{|t| \ge y^2/\log^2 y} \frac{1}{\lam t^2} dt 
\ll \frac{{\Bbb E}(L(1,X;y)^{k})}{(e^{\gamma} \tau)^k} \frac{\log^2 y}
{\lam y^2}.
$$
Using Lemma 3.2 and as $(e^{\lam s}-1)/(\lam s)\ll 1$ and 
$(e^{\lam s}-e^{-\lam s})/s \ll \lam$ 
we may see easily that the contribution of the 
initial segment of the integral is 
$\ll {\Bbb E}(L(1,X;y)^{k})(e^{\gamma} \tau)^{-k} \lam \sqrt{k\log k}$.  
Thus we conclude that 
$$
{\Phi}(\tau e^{-\lam};y) -{\Phi}(\tau;y) 
\ll \frac{{\Bbb E}(L(1,X;y)^{k})}{(e^{\gamma} \tau)^k} 
\Big( \frac{\log^2 y}{\lam  y^2} + \lam \sqrt{k\log k}\Big). \tag{3.8}
$$

Consider next the integral in (3.6).  We split the 
integral into three parts: when $|t|\le \sqrt{k} (\log k)^2$, 
when $\sqrt{k} (\log k)^2\le |t|\le y^2/\log^2 y$ and when $|t|>y^2/\log^2 y$. 
In the third case the integrand is $\ll 
{\Bbb E}(L(1,X;y)^k)(e^{\gamma}\tau)^{-k}/(\lam |t|^2)$, and 
so the contribution of this segment of the integral is 
$\ll {\Bbb E}(L(1,X;y)^k)(e^{\gamma}\tau)^{-k}
\log^2 y/(\lam y^2)$. Using Lemma 3.2 we see that 
the contribution of the second segment of the integral is 
bounded by ${\Bbb E}(L(1,X;y)^k)(e^{\gamma}\tau)^{-k} \exp(-c \log^3 k)/\lam$.
Lastly, using (3.5) and our definition of $k$, we get that the initial 
segment of the integral contributes 
$$
\align
&\frac{{\Bbb E}(L(1,X;y)^k)}{(e^{\gamma }\tau)^k} \frac{e^{\lam k}-1}
{\lam k^2}\Big( 1+ O\Big(\frac{\log^5 k}{\sqrt{k}}\Big)\Big) 
\frac{1}{2\pi} \int_{|t|\le \sqrt{k} \log^2 k} \exp\Big(-\frac{t^2}{2}R(k;y)
\Big)dt\\
&= \frac{{\Bbb E}(L(1,X;y)^k)}{(e^{\gamma }\tau)^k} \frac{e^{\lam k}-1}
{\lam k^2}\frac{1}{\sqrt{2\pi R(k;y)}}
\Big( 1+ O\Big(\frac{\log^5 k}{\sqrt{k}}\Big)\Big).
\\
\endalign
$$
Observe that $R(k;y) \asymp 1/(k\log k)$ by (3.2), and so 
choosing $\lam= k^{-\frac 54}$ we obtain from the above that 
$$
{\Phi}(\tau e^{-\lam};y) \ge 
\frac{{\Bbb E}(L(1,X;y)^k)}{(e^{\gamma }\tau)^k} \frac{1}{k\sqrt{2\pi R(k;y)}}
\Big( 1+ O\Big(\frac{\log^2 k}{k^{\frac 14}}\Big)\Big) 
\ge {\Phi}(\tau;y).
$$
Further with this same choice of $\lam$ we obtain by (3.8) that 
$$
{\Phi}(\tau e^{-\lam};y)-{\Phi}(\tau;y) \ll 
\frac{{\Bbb E}(L(1,X;y)^k)}{(e^{\gamma }\tau)^k} \frac{1}{k\sqrt{2\pi R(k;y)}}
\frac{\log^2 k}{k^{\frac 14}},
$$
and so the first part of the Theorem follows.

Now noting again that $R(k;y) \asymp 1/(k\log k)$, and using the 
first part of the Theorem, we may write (3.8) as 
$$
{\Phi}(\tau e^{-\lam};y) -{\Phi}(\tau;y) \ll {\Phi}
(\tau;y) \Big( \frac{\sqrt{k}\log^2 y}{y^2 \lam}+ \lam k \Big), 
$$ 
which gives the second part of the Theorem when $\log y/(y e^{\frac 14 \tau})
\asymp \log y/(y k^{\frac 14})\le \lam  \le e^{-\tau}$.  Since 
$\Phi(\tau e^{-\lam};y)-\Phi(\tau;y)$ is a non-decreasing 
function of $\lam$, thus for $\lam \le \log y/(y e^{\frac 14 \tau})$ we have 
$\Phi(\tau e^{-\lam};y)-\Phi(\tau;y)\ll \Phi(\tau;y) e^{3\tau/4}\log y / y$
as desired.

\enddemo

In order to extract simpler asymptotics for $\Phi(\tau;y)$ we 
now need some understanding of the asymptotic nature of 
${\Bbb E}(L(1,X;y)^k)$.  
If $k\le y$ then using (3.4) for $p\le k/\log k$, and (3.3) 
for $k/\log k\le p\le y$ we deduce that 
$$
{\Bbb E}(L(1,X;y)^k) \asymp \prod_{p\le k/\log k} 
\frac{p}{2(p+1)} \Big( 1-\frac{1}{p}\Big)^{-k} \prod_{k/\log k\le p\le y} 
\cosh (k/p). \tag{3.9}  
$$
Now note that 
$$
\prod_{k/\log k \le p\le k} e^{-k/p} \Big( 1-\frac1p\Big)^{-k}  
\asymp \prod_{k/\log k \le p\le k} e^{O(k/p^2)} \asymp 1,
$$ 
and also that by partial summation using the prime number theorem 
$$
\prod_{k/\log k \le p\le k} e^{-k/p} \cosh(k/p) \prod_{y\ge p>k} \cosh(k/p) 
= \exp\Big( \frac{k}{\log k} \Big(\int_{k/y}^{\infty} g_1(t)\frac{dt}{t^2} 
+ O\Big( \frac{1}{\log k}\Big)\Big)\Big), 
$$
where $g_1(t) =\log \cosh (t) -t$ if $t>1$ and $\log \cosh (t)$ if $t\le 1$.
Integration by parts shows easily that 
$$
\int_{k/y}^{\infty} g_1(t)\frac{dt}{t^2} = C_1(k/y) -1+\frac{g_1(k/y)}{k/y}. 
$$
Using these estimates in (3.9) we conclude that 
$$
 {\Bbb E}(L(1,X;y)^k) = \prod_{p\le k}\Big(1-\frac 1p\Big)^{-k} 
\exp\Big( \frac{k}{\log k} \Big(C_1(k/y) -1+ \frac{\log \cosh(k/y)}{k/y} 
+O\Big(\frac{1}{\log k}\Big)\Big)\Big). \tag{3.10}
$$
Similarly we find that 
$$ 
{\Bbb E}(L(1,X;y)^{-k})= \prod_{p\le k}\Big(1+\frac 1p\Big)^{k} 
\exp\Big( \frac{k}{\log k} \Big(C_1(k/y) -1+ \frac{\log \cosh(k/y)}{k/y} 
+O\Big(\frac{1}{\log k}\Big)\Big)\Big). 
$$
From these estimates and 
Theorem 3.1 we may deduce the following Corollary which contains 
Proposition 1.  

\proclaim{Corollary 3.3}  
Let $y$ be large, and let $\tau$ be large with $\tau \le \log y - 1$.
Then 
$$
{\Phi}(\tau;y) = \exp\Big( - \frac{e^{\tau-C_1}}{\tau}
\Big(1+O\Big(\frac {e^{\tau}}{y} +\frac{1}{\tau}\Big)\Big)\Big).
$$
The same asymptotic holds for ${\Psi}(\tau;y)$.  
\endproclaim

\demo{Proof}  We take $k=k_{\tau,y}$ as in Theorem 3.1.  By Mertens' 
theorem, (3.10), and Theorem 3.1 we see that 
$$
{\Phi}(\tau) = \Big( \frac{\log k}{\tau}\Big)^k \exp\Big( 
\frac{k}{\log k} \Big(C_1(k/y) -1+ \frac{\log \cosh(k/y)}{k/y} 
+O\Big(\frac{1}{\log k}\Big)\Big)\Big),
$$
and using (3.1) the definition of $k$ this is 
$$
= \exp\Big( 
\frac{k}{\log k} \Big( -1+ \frac{\log \cosh(k/y)}{k/y} 
+O\Big(\frac{1}{\log k}\Big)\Big)\Big).
$$
Now $y \ge k\asymp e^{\tau}$ and so (3.1) gives that 
$I(k;y) = \log_2 k +\gamma + C_1(k/y)/\tau + O(1/\tau^2) 
= \log_2 k + \gamma + C_1/\tau + O(e^{\tau}/(\tau y)+1/\tau^2)$ 
from which it follows that $k=e^{\tau -C_1} (1+O(e^{\tau}/y+1/\tau))$.  
Using this above we obtain the asymptotic for $\Phi(\tau;y)$, 
and the proof for $\Psi(\tau;y)$ is similar.

\enddemo

\head 4.  Estimates for real character sums \endhead

\noindent In this section we collect together some estimates 
for $\sumf_{|d|\le x} \L{d}{n}$ for individual $n$, and 
also on average over $n$ in a dyadic interval.  
We begin with a simple application of the P{\' o}lya-Vinogradov 
inequality. 

\proclaim{Lemma 4.1}  If $n$ is a positive integer, not a
perfect square, then 
$$
\Big|\sumf_{|d| \le x} \L{d}{n} \Big| \ll x^{\frac 12} n^{\frac 14} 
(\log n)^{\frac 12}.
$$
\endproclaim
\demo{Proof}  We shall confine our attention to $d$ positive
and $\equiv 1\pmod 4$:  the cases $d$ negative, or 
$d\equiv 8$, or $12 \pmod{16}$ are handled similarly.  
Thus we seek to bound 
$$
\sum\Sb d\le x \\ d\equiv 1 \pmod 4 \endSb \mu^2(d) \L{d}{n} 
= \frac{1}{2} \sum_{\psi \pmod 4} \sum\Sb d\le x\endSb \psi(d) \mu^2(d) 
\L{d}{n}.
$$
Writing $\mu^2(d)=\sum_{l^2 |d} \mu(l)$, and using 
the P{\' o}lya-Vinogradov inequality (which is applicable 
since $\fracwithdelims() {\cdot}{n} \psi(\cdot)$ is 
a non-principal character of conductor at most $4n$) 
we get that the above is 
$$
\align
&\ll  \sum_{\psi \pmod 4} \sum_{l \le \sqrt{x}} 
\Big| \sum\Sb d\le x \\ l^2 |d \endSb \L{d}{n} \psi(d)\Big| 
\ll \sum_{\psi \pmod 4}  
\sum_{l\le \sqrt{x}} \Big|\sum_{m\le x/l^2} 
\fracwithdelims() mn \psi(m)\Big|\\
&\ll \sum_{l\le \sqrt{x}} \ \min \Big( \frac{x}{l^2},\sqrt{n} \log n\Big) 
\ll  x^{\frac 12} n^{\frac 14} (\log n)^{\frac 12}, 
\\
\endalign
$$
which proves the Lemma.
\enddemo

Since $\sumf_{|d|\le x} \L{d}{n}$ is trivially $\ll x$ we 
see that Lemma 4.1 furnishes a non-trivial bound 
only when $x \ge \sqrt{n} \log n$.  While this range 
can be improved using Burgess' character sum estimates, in general 
non-trivial bounds are known only when $x$ is larger than some 
fixed power of $n$.  In the special case that $n$ is 
smooth (that is composed only of small prime factors) then 
we may use a remarkable result of Graham and Ringrose 
to obtain non-trivial bounds in the range $x\ge n^{\epsilon}$.  This 
will be a crucial ingredient in Section 6 below.

Suppose $\chi$ is a non-principal character 
$\pmod q$ where $q/(4,q)$ is square-free, and suppose 
$p$ is the largest prime factor of $q$. Then
for any integer $l \ge 2$ and with $L=2^l$ we have
$$
\sum_{n\leq N} \chi(n) \ll N^{1-\frac{l}{8L}}  
p^{\frac 13} q^{\frac1{7L}} d(q)^{\frac{l^2}{L}}.\tag{4.1}
$$ 
This is a consequence of Theorem 5 of Graham and Ringrose [8], 
using there that $\log q\ll p$ and 
$\prod_{p|q} (1+1/p)\ll \log\log q$, and simplifying their estimates.

\proclaim{Lemma 4.2} Let $n$ be a positive integer not a perfect square.  
Write $n=n_0 \square$ where $n_0>1$ is square-free, and suppose 
all prime factors of $n_0$ are below $P$.  Let $l\ge 2$ be an integer 
and put $L=2^l$.  Then 
$$
\sumf_{|d|\le x} \L{d}{n} \ll 
 x^{1-\frac{l}{8L}} \prod_{p|n} \Big(1+\frac{1}{p^{1-l/8L}}\Big) 
P^{1/3} n_0^{\frac{1}{7L}} d(n_0)^{\frac{l^2}{L}}.
$$
\endproclaim 

\demo{Proof}  We shall confine our attention to $d$ positive
and $\equiv 1\pmod 4$:  the cases $d$ negative, or $d \equiv 8$, or 
$12 \pmod {16}$ are handled similarly.  Thus we seek to bound
$$
\sum\Sb d\le x\\ d\equiv 1\pmod 4\endSb \mu^2(d) \L{d}{n} 
= \sum\Sb d\le x\\ d\equiv 1\pmod 4\\ (d,n)=1\endSb \mu^2(d) \L{d}{n_0} 
= \frac 12 \sum_{\psi \pmod 4} \sum\Sb d\le x\\(d,n)=1\endSb 
\psi(d) \L{d}{n_0} \mu^2(d).
$$
Note that $\sum_{a^2 |d, (a,n)=1} \mu(a) \sum_{b|(d,n)} \mu(b) =1$ 
if $d$ is square-free and coprime to $n$, and $0$ otherwise.
Hence the above is 
$$
\align
&\ll \sum_{\psi \pmod 4} \sum_{b|n } \mu^2(b) 
\sum\Sb a^2 b\le x\\ (a,n)=1\endSb  \mu^2 (a)
\Big| \sum\Sb d\le x\\ a^2 b|d \endSb \L{d}{n_0} \psi(d)\Big|
\\
&\ll \sum_{\psi \pmod 4} \sum_{b|n } \mu^2(b) 
\sum\Sb a^2 b\le x\\ (a,n)=1\endSb  \mu^2 (a)
\Big| \sum_{m\le x/a^2b} \psi(m)\L{m}{n_0}\Big|\\
&\ll  \sum_{b|n } \mu^2(b) 
\sum\Sb a^2 b\le x\\ (a,n)=1\endSb  \mu^2 (a)
 \Big(\frac{x}{a^2 b}\Big)^{1-\frac{l}{8L}}
P^{1/3} n_0^{\frac{1}{7L}} d(n_0)^{\frac{l^2}{L}} ,\\
\endalign
$$
by (4.1), since $\psi(\cdot) \L{\cdot}{n_0}$ is a non-principal character 
$\pmod {n_0}$ or $\pmod {4n_0}$, which yields the Lemma.
\enddemo

We now give results bounding $\sumf_{|d|\le x}\L{d}{n}$ on average over 
$n$.  To do this we shall use the following 
consequence of a simple large sieve estimate for real characters.  

\proclaim{Lemma 4.3}  For non-zero integers $m \equiv 0,\  1\pmod 4$, 
and natural numbers $n$, let $a_m$ and 
$b_n$ denote arbitrary complex numbers, and set $a_0 =0$ and 
$a_m=0$ if $m \equiv 2, 
\ 3 \pmod 4$.  
Then 
$$
\sum_{n \le N} \Big|\sum\Sb |m|\le M\endSb 
 a_m \L{m}{n}\Big|^2 \ll 
N\sum\Sb |m_1|, |m_2| \le M\\ m_1m_2=\square \endSb |a_{m_1} a_{m_2}| 
+ M\log M \Big(\sum_{|m|\le M} |a_m|\Big)^2, \tag{4.2}
$$
and 
$$
\sum\Sb |m|\le M\endSb 
 \Big| \sum_{n\le N} b_n \L{m}{n}\Big|^2 
\ll M \sum\Sb n_1, n_2 \le N\\ n_1n_2 =\square \endSb |b_{n_1}b_{n_2}| 
+ N\log N \Big( \sum_{n\le N} |b_n|\Big)^2. \tag{4.3}
$$
Alternative bounds are
$$
\sum_{n \le N} \Big|\sum\Sb |m|\le M\endSb 
 a_m \L{m}{n}\Big|^2 
\ll (M^2 \sqrt{N} + N^{\frac 32} M) \log (MN) \Big( 
\sum_{|m|\le M} \frac{d(|m|) |a_m|^4 }{|m|} \Big)^{\frac 12}, \tag{4.4}
$$
and 
$$
\sum\Sb |m|\le M\endSb 
 \Big|\sum_{n\le N} b_n \L{m}{n}\Big|^2 \ll 
(N^2 \sqrt{M}+M^{\frac 32} N) 
\log (MN) \Big(\sum_{n\le N} \frac{d(n)|b_n|^4}{n}\Big)^{\frac 12}.
\tag{4.5}
$$
\endproclaim
\demo{Proof}  The first two estimates are simple consequences 
of the P{\'o}lya-Vinogradov inequality.  Taking (4.2)  
for instance, the desired quantity is 
$$
\sum_{|m_1|, |m_2| \le M} a_{m_1} \overline{a_{m_2}} \sum_{n\le N} 
\L{m_1 m_2}{n}.
$$
Now $m_1 m_2 \equiv 0, \ 1 \pmod 4$ so that $\L{m_1 m_2}{\cdot}$ is a 
character $\pmod{|m_1m_2|}$, and this 
character is principal when $m_1m_2 =\square$ giving rise to the 
first term in the RHS of (4.2), 
and the character is non-principal when $m_1 m_2 \neq \square$ giving 
rise (via  P{\'o}lya-Vinogradov) by the second term there.  The proof 
of (4.3) is similar.

We now prove the second estimates: again we give only the 
case (4.4), the proof of (4.5) being entirely 
similar.  Note that our 
desired expression is 
$$
\sum_{|m| \le M^2} \Big( \sum\Sb m_1 m_2 =m\\ |m_1|, |m_2| \le M\endSb 
a_{m_1} \overline{a_{m_2}} \Big) \sum_{n\le N} \L{m}{n}
$$
which is by Cauchy-Schwarz 
$$
\le \Big( \sum_{|m|\le M^2} \Big|  \sum\Sb m_1 m_2 =n\\ |m_1|, |m_2| 
\le M\endSb 
a_{m_1} \overline{a_{m_2}} \Big|^2 \Big)^{\frac 12} 
\Big( \sum_{|m|\le M^2} \Big|\sum_{n\le N} \L{m}{n}\Big|^2 \Big)^{\frac 12}.
$$
Applying (4.3) 
we see that 
the second factor above is $\ll (M^2 N \log N + N^3 \log N)^{\frac 12}$.  
Further as 
$$
\Big|  \sum\Sb m_1 m_2 =m\\ |m_1|, |m_2| \le M\endSb 
a_{m_1} \overline{a_{m_2}} \Big|^2 \le \Big|\sum\Sb d |m 
\\ |d|\le M\endSb 
|a_d|^2 \Big|^2 \le 2 d(|m|) \sum\Sb d|m \\ |d|\le M\endSb |a_d|^4,
$$
by Cauchy-Schwarz, we get that the first factor above is 
$$
\ll \Big(\sum_{|d|\le M} |a_d|^4 \sum\Sb |m|\le M^2 \\ d|m\endSb d(|m|)
\Big)^{\frac 12}   
\ll \Big( M^2 \log M \sum_{|d|\le M} \frac{|a_d|^4 d(|d|)}{|d|} 
\Big)^{\frac 12},
$$
completing the proof of the Lemma. 
 
\enddemo

Applying the above lemma we obtain the following estimate 
for the $2k$-th moment of $\sumf_{|d|\le x} \L{d}{n}$ averaged over 
$n$.  This will form the key input in our first approach to 
the moments of $L(1,\chi_d)$; see Section 5 below.  

\proclaim{Lemma 4.4}  Uniformly for all integers $k\ge 1$ we have 
$$
 \sum_{n\le N} \Big| \sumf_{|d|\le x} \L{d}{n}\Big|^{2k} 
\ll (x^{2k} N^{\frac{1}{2}} + x^{k} N^{\frac{3}{2}} ) 
(2k\log x)^{2k^5}, 
$$
and also 
$$
\ll (x^{k}N +x^{3k} ) (2k\log x)^{3k^2}.
$$
\endproclaim
\demo{Proof}  Write 
$$
\Big( \sumf_{|d|\le x} \L{d}{n}\Big)^k = \sum_{|m| \le x^k } 
a_m \L{m}{n},
$$ 
where $a_m=0$ unless $1\le | m|\le x^k$ is $\equiv 0, \ 1\pmod 4$, and 
$0 \le a_m \le d_k(|m|)$ for such $m$.  Note also 
that $\sum_{|m|\le x^k} a_m = (\sumf_{|d|\le x} 1)^k \le x^k$.  

Then from (4.2) we get that 
$$
\sum_{n\le N} \Big| \sumf_{|d|\le x} \L{d}{n}\Big|^{2k} 
\ll N \sum\Sb m_1, m_2 \le x^k\\ m_1 m_2 =\square \endSb 
 d_k(m_1)d_k(m_2) + x^k \log (x^k) \Big(\sum_{|m| \le x^k }a_m\Big)^2.
\tag{4.6}
$$
The second term above is $\le x^{3k} \log (x^k)$.  As 
for the first term above, note that $m_1 m_2 =\square$ means that 
we may write $m_1= n\alpha^2$, and $m_2 =n\beta^2$ where 
$n\le x^k$ and $\alpha$ and $\beta$ are $\le \sqrt{x^k/n}$.  
Since $d_k(n\alpha^2) d_k(n\beta^2) \le d_k(n)^2 d_k(\alpha^2) d_k(\beta^2) 
\le d_{k^2}(n) d_{k^2}(\alpha) d_{k^2}(\beta)$, and $n\alpha\beta \le x^k$, 
we see that 
$$
\align
\sum\Sb m_1, m_2 \le x^k\\ m_1 m_2 =\square \endSb 
 d_k(m_1)d_k(m_2) &\le x^{k} \sum_{n, \alpha, \beta \le x^k} \frac{d_{k^2} (n) 
d_{k^2} (\alpha) d_{k^2}(\beta)}{n\alpha \beta} 
\le x^k \Big( \sum_{n\le x^k } \frac{d_{k^2} (n)}{n}\Big)^3 \\
&\le x^k \Big(\sum_{n\le x^k} \frac{1}{n}  \Big)^{3k^2}
\le x^k (2k \log x)^{3k^2}.
\\
\endalign
$$
Inputing these estimates in (4.6) we get 
the second bound claimed in the Lemma.  

The first bound of the Lemma is trivial when $N\ge x^{2k}$ since the sum is
$\leq Nx^{2k}\leq x^kN^{3/2}$. Suppose now that $N\le x^{2k}$.  
Using (4.4) we get that 
$$
\sum_{n\le N} \Big| \sumf_{|d|\le x} \L{d}{n}\Big|^{2k} 
\ll (x^{2k} \sqrt{N} + N^{\frac 32} x^k) (3k\log x) 
\Big( \sum_{n\le x^k} \frac{d(n)d_k(n)^4}{n}\Big)^{\frac 12}.
$$
Since $d(n)d_k(n)^4\leq d(n)d_{k^4}(n)\leq d_{k^4+1}(n) \leq d_{2k^5}(n)$
we have
$$
\sum_{n\le x^k} \frac{d(n)d_k(n)^4}{n}
\le \sum_{n\le x^k} \frac{d_{2k^5}(n)}{n}\le  
\Big(\sum_{n\le x^k} \frac 1n\Big)^{2k^5} \le (2k\log x)^{2k^5},
$$
and the Lemma follows.

\enddemo

Lastly we give a bound for $\sumf_{|d|\le x} \L{d}{n}$ for non-square 
integers $n$, conditional on the Generalized Riemann Hypothesis.  

\proclaim{Lemma 4.5}  Assume GRH. Fix $\varepsilon>0$.
Let $\chi$ be a non-principal character $\pmod q$.  
Then for $x\ge 2$
$$
\sum_{m\le x} \chi(m) \ll x^{\frac 34+\varepsilon} \exp( (\log q)^{1/2-\varepsilon}).
$$
Consequently, for any positive integer $n$ which is not a perfect 
square we have
$$
\sumf_{|d|\le x} \L{d}{n} \ll x^{\frac 34+\varepsilon} \exp((\log n)^{1/2-\varepsilon}).
$$ 
\endproclaim 

\demo{Proof}  We may suppose that $x\le q$, and also that 
$x$ is half an odd integer.  Then by Perron's formula 
$$
\align
\sum_{n\le x} \chi(n) &= \frac{1}{2\pi i} \int_{1+\frac{1}{\log x} 
-i\infty}^{1+\frac{1}{\log x} + i\infty} L(s,\chi ) x^s \frac{ds}{s} \\
&= \frac{1}{2\pi i} \int_{1+\frac{1}{\log x} 
-iq}^{1+\frac{1}{\log x} + iq} L(s,\chi ) x^s \frac{ds}{s} 
+O\Big(\frac{x\log^2 x}{q}\Big).
\\
\endalign
$$
Appealing to Lemma 2.1 with $\sigma_0=\frac 12$ there (by GRH) and 
$y=(\log q)^2(\log_2 q)^6$, 
we deduce that if Re$(s) \ge \frac 34 +\varepsilon$, and Im$(s) \le q$ then 
$\log L(s,\chi) \leq (1/2) (\log q)^{1/2-\varepsilon}$.  
Using this bound, and moving the line of integration above to the 
Re$(s)=\frac 34+\varepsilon$ line we obtain the first part of the Lemma.  
To obtain the second part of the lemma, we proceed along 
the lines of the proof of Lemma 4.1, replacing the 
use of P{\' o}lya-Vinogradov by the GRH bound 
for character sums above.

\enddemo 

Different bounds may be found by moving to 
the line Re$(s)=\sigma$ for any $1>\sigma >\frac 12$,  but the present 
bound is adequate for our applications.

\head 5.  Moments of $L(1,\chi_d)$:  Proof of Theorem 2 \endhead

\noindent By Lemma 2.3 we see that with $Z= \exp( (\log x)^{10})$ 
$$ 
\sumf\Sb |d|\le x \\ d\notin {\Cal L}\endSb L(1,\chi_d)^z 
= \sum_{n=1}^{\infty} \frac{d_z(n)}{n} e^{-n/Z} \sumf\Sb |d|\le x \\ 
d\notin {\Cal L}\endSb \L{d}{n} + O(\log x).
$$
Since there are $\ll \log x$ discriminants $d\in {\Cal L}$, $|d|\le x$, we 
see that the above is 
$$
\align
&\sum_{n=1}^{\infty} \frac{d_z(n)}{n} e^{-n/Z} \sumf\Sb |d|\le x \endSb 
 \L{d}{n} + O\Big( \log x \sum_{n=1}^{\infty} \frac{|d_z(n)|}{n} e^{-n/Z}\Big) 
\\
= &\sum_{n=1}^{\infty} \frac{d_z(n)}{n} e^{-n/Z} \sumf\Sb |d|\le x \endSb 
 \L{d}{n} + O(x^\epsilon), \tag{5.1}
\\
\endalign
$$
using $|d_z(n)|\le d_{|z|}(n)$ and (2.4), 
and since $|z|\le \log x/(\log_2 x)^2$.  

We now handle the contribution of the terms $n=\square$ to (5.1), 
which gives the main term.  When $n=m^2$ we have
$$
\sumf_{|d|\le x} \L{d}{m^2} =\sumf\Sb |d|\le x\\ (d,m)=1\endSb 
1 = \frac{6}{\pi^2}x \prod_{p|m} \L{p}{p+1} + O(x^{\frac 12+\epsilon} 
d(m)).
$$ 
Thus the contribution of such terms to (5.1) is
$$
\frac{6}{\pi^2}x \sum_{m=1}^{\infty} \frac{d_z(m^2)}{m^2} \prod_{p|m} 
\L{p}{p+1} e^{-m^2/Z} + O\Big(x^{\frac 12+\epsilon} \sum_{m=1}^{\infty} 
\frac{|d_z(m^2) d(m)|}{m^2}  e^{-m^2/Z} \Big).
$$
Since $|d_z(m^2)d(m)| \le d_{2|z|+2}(m)^2$ (see section 2) 
and $ e^{-m^2/Z} \leq 1$ the error term above is 
$$
\align
&\ll x^{\frac 12+ \epsilon} \sum_{m=1}^{\infty} 
\frac{d_{2|z|+2}(m)^2}{m^2}  
\ll x^{\frac 12+\epsilon} \prod_{p} \Big( \int_0^1 \Big|1- \frac{e(\theta)}{p}
\Big|^{-4(|z|+1)} d\theta \Big)
\\
&\ll 
x^{\frac 12+\epsilon} \prod_{p\le |z|+2} \Big( 1-\frac 1p\Big)^{-4(|z|+1)} 
\prod_{p >|z|+2} \Big(1+O\Big(\frac{|z|^2}{p^2}\Big)\Big)
\\
&\ll x^{\frac 12+\epsilon} (10\log (|z|+2))^{4|z|+4} 
\ll x^{\frac 23},
\\
\endalign
$$
say, since $|z| \le \log x/(\log_2 x)^2$.  
Further  since $1-e^{-t} \leq t^{1/4}$ for any $t>0$, and
$|d_z(m^2) |\le d_{(|z|+1)^2} (m)$, we have
$$
\sum_{m=1}^{\infty} \frac{|d_z(m^2)| }{m^2} (1-e^{-m^2/Z}) 
\le  \sum_{m=1}^{\infty} \frac{d_{(|z|+1)^2}(m)}{m^2}
\left( \frac{m^2}{Z}\right)^{1/4}
=  \frac{\zeta(3/2)^{(1+|z|)^2}}{Z^{1/4}} \leq \frac 1x .
$$
We conclude that the contribution of the terms $n=\square$ to (5.1) 
is 
$$
\frac{6}{\pi^2} \sum_{m=1}^{\infty} \frac{d_z(m^2)}{m^2} \prod_{p|m} 
\L{p}{p+1} + O(x^{\frac 23}) 
= \frac 6{\pi^2} {\Bbb E} (L(1,X)^z) + O(x^{\frac 23}). 
\tag{5.2}
$$

We now handle the $n\neq \square$ terms in (5.1).  Let 
$k$ be a positive integer $\geq 2$ to be fixed shortly.  Using H{\" o}lder's 
inequality with exponents $\alpha:= 2k$ and 
$\beta:=2k/(2k-1)$ (so that $1/\alpha + 1/\beta =1$) 
we see that 
$$
\sum\Sb n=1\\ n\neq \square\endSb^{\infty} 
\frac{|d_z(n)|}{n} e^{-\frac{n}{Z}} 
\Big| \sumf_{|d| \le x} \fracwithdelims() dn \Big| 
\leq
\Big( 
\sum_{n=1}^{\infty} \frac{|d_z(n)|^\beta}{n} e^{-\frac{n}{Z}} 
\Big)^{\frac{1}\beta} 
\Big( \sum\Sb n=1 \\ n\neq \square \endSb^{\infty} 
\frac{e^{-\frac{n}{Z}}}{n} \Big| \sumf_{|d| \le x}
 \fracwithdelims() dn 
\Big|^{2k} \Big)^{\frac{1}{2k}}.
\tag{5.3}
$$
Since $|d_z(n)|^\beta \le d_{|z|^\beta} (n)$,  (2.4) gives 
that the first factor here is 
$$
\leq (\log 3Z)^{(|z|+1)^\beta} 
= \exp\left( 20(|z|+1)^\beta \log \log x\right). \tag{5.4} 
$$
We split the second factor into dyadic blocks:
$$
\sum\Sb n=1 \\ n\neq \square \endSb^{\infty} 
\frac{e^{-\frac{n}{Z}}}{n} 
\Big| \sumf_{|d| \le x} \fracwithdelims() dn 
\Big|^{2k} 
\leq \sum_{j=0}^{\infty} \frac{e^{-2^{j}/Z}}{2^{j}} \sum\Sb n=2^{j} \\ 
n \neq \square\endSb^{2^{j+1}-1} \Big|\sumf_{|d| \le x}
 \fracwithdelims() dn \Big|^{2k}.
$$
To bound the contribution of the dyadic block $[2^j,2^{j+1})$, 
we use the bound of Lemma 4.1 when $2^{j} \le x^{2k/(k+1)}$,
the first bound of Lemma 4.4 when $x^{2k/(k+1) } < 2^j \le x^{4k/3}$
and the second bound of Lemma 4.4 when $2^j \ge x^{4k/3}$.  
Summing these bounds we deduce that the above is 
$\ll x^{k(2k+1)/(k+1)} (2k\log x)^{2k^5}$.  
Thus the second factor in (5.3) is 
$\ll x^{(2k+1)/(2k+2)}(2k\log x)^{k^4}$.   
Using this with (5.4) we 
conclude that the contribution of the $n\neq \square$ 
terms is 
$$
\ll x^{1-\frac{1}{2k+2}} (2k\log x)^{k^4}
 \exp\left( 20(|z|+1)^\beta \log \log x\right). 
$$ 
Choose $k=[\log_2 x]$ and recall that $|z|\le \log x/(500 (\log_2 x)^2)$. 
Then the above error is $\ll x \exp(-\log x/(5\log_2 x))$, and this 
when combined with (5.2) proves our Theorem.

\head 6. Moments of short Euler products:  Proof of Theorem 3
\endhead

\noindent 

\proclaim{Theorem 6.1}  Let $2\le y \le \exp(\sqrt{\log x})$, 
and let $z$ be a complex 
number.   Write  
$$
\sumf_{|d|\le x} L(1,\chi_d;y)^z = \frac{6}{\pi^2} x {\Bbb E}(L(1,X;y)^z) 
+ \theta x {\Bbb E}(L(1,X;y)^{\text{Re }z}) E(z,y), 
$$
where $|\theta| =1$, and $E(z,y)$ is a positive real number. 
If $y \ge 4|z|+4$ then 
$$
E(z,y) \ll y^{\frac 13} \exp\Big( -\frac{\log x \log_2 y}{12\log y} + 
\frac{30|z|}{\log (4|z|+4)} + 10|z| \log \Big(\frac{\log y}{\log (4|z|+4)}
\Big)\Big),
$$
and if $y\le 4|z|+4$ then 
$$
E(z,y) \ll x^{-\frac 1{40}}, \qquad \text{if   } y \le \frac{\log x}{2}, 
$$
and 
$$
E(z,y) \ll \exp\Big( - \frac{\log x}{(\log_2 x)^{\frac 34}} \Big), 
\qquad \text{if  } y\le \frac 13 \log x \log_3 x.
$$
\endproclaim

\demo{Proof}  Let $\delta = \pm 1$ indicate the sign of $\text{Re }z$.
Denote $\sumf_{|d|\le x} \chi_d(n)$ by $S(x;n)$, and define
${\Cal S}(y)$ to be the set of integers whose prime factors are all $\leq y$.
  Then 
$$
\sumf_{|d|\le x} L(1,\chi_d;y)^z = \sumf_{|d|\le x} 
\sum\Sb n=1 \\ n\in {\Cal S}(y)\endSb^{\infty} \frac{d_z(n)}{n} \chi_d(n) 
= \sum\Sb n=1 \\ n\in {\Cal S}(y)\endSb^{\infty} \frac{d_z(n)}{n}
S(x;n). 
$$
Decompose $n=n_1 n_2^2 n_3^2$ where $n_1$ and $n_2$ are 
square-free, with $(n_1,n_2)=1$, and $p|n_3 \implies p | n_1n_2$: 
that is $n_1$ is the product of all primes dividing $n$ to an 
odd power, and $n_2$ is the product of all primes dividing $n$ to an 
even power $\ge 2$. Observe that $S(x;n_1n_2^2n_3^2) = S(x;n_1 n_2^2)$.  
Thus our sum above is 
$$
\align
&\sum\Sb n_1 \in {\Cal S}(y) \endSb \mu^2(n_1) 
\sum\Sb n_2 \in {\Cal S}(y) \\ (n_1,n_2)=1 \endSb 
\mu^2(n_2) S(x;n_1 n_2^2) \sum\Sb p|n_3 \implies p|n_1n_2\endSb 
\frac{d_z(n_1n_2^2 n_3^2)}{n_1n_2^2 n_3^2} 
\\
= &\sum\Sb n_1 \in {\Cal S}(y) \endSb \mu^2(n_1) 
\sum\Sb n_2 \in {\Cal S}(y) \\ (n_1,n_2)=1 \endSb 
\mu^2(n_2) 
S(x;n_1 n_2^2) \prod_{p|n_1} s_p(z) \prod_{p|n_2} (c_p(z)-1), \tag{6.1}
\\
\endalign
$$ 
where 
$$
s_p(z) := \frac{(1-1/p)^{-z}-(1+1/p)^{-z}}{2}, \qquad 
\text{and} \qquad c_p(z):= \frac{(1-1/p)^{-z}+(1+1/p)^{-z}}{2}.
$$
Remark that $|s_p(z)|$ and $|c_p(z)|$ are always 
$\le (1-\delta/p)^{-\text{Re }z}$, and that 
when $p > 4|z|+4$ we have $|s_p(z)| \le 2|z|/p$, and $|c_p(z)-1| 
\le 2|z|^2/p^2$.  To evaluate (6.1), we now 
distinguish the cases $n_1=1$ which give rise to 
the main term, and $n_1>1$ which contribute the error term. 
 
We first handle the terms $n_1 >1$.  By Lemma 4.2 we 
see that for any integer $l\ge 2$ and with $L=2^l$ we
have 
$$
|S(x;n_1n_2^2)| \ll x^{1-\frac{l}{8L}} \prod_{p|n_1n_2} \Big(1 +
\frac{1}{p^{1-l/8L}}\Big) y^{\frac 13} n_1^{\frac{1}{7L}} 
d(n_1)^{\frac{l^2}{L}}. 
$$
Using this estimate in (6.1) we see that 
the contribution of the terms $n_1>1$ is 
$$
\align
&\ll x^{1-\frac{l}{8L}} y^{\frac 13} \prod_{p\le y} 
\Big( 1+ 2^{\frac{l^2}{L}} p^{\frac 1{7L}} \Big(1+\frac{1}{p^{1-1/8L}}
\Big) |s_p(z)| +  \Big(1+\frac{1}{p^{1-1/8L}}
\Big)|c_p(z)-1|\Big)
\\
&\ll  x^{1-\frac{l}{8L}} y^{\frac 13} \prod_{p \le \min(4|z|+4,y)} 
\Big( 8p^{\frac 1{7L}} \Big(1-\frac {\delta}p\Big)^{-\text{Re } z} \Big) 
\prod_{\min(4|z|+4, y)\le p\le y} \Big(1+\frac{8|z|p^{1/7L}}{p}\Big).
\\
\endalign
$$

We now handle the terms arising from $n_1=1$.  Note that 
$S(x;n_2^2)$ counts the number of fundamental discriminants $|d|\le x$ 
that are coprime to $n_2$.  Thus we see 
easily that $S(x;n_2^2) = \frac{6}{\pi^2 }x \prod_{p|n_2} \L{p}{p+1} 
+ O(\sqrt{x} d(n_2))$.  Hence the contribution of the terms $n_1=1$ to 
(6.1) equals 
$$
\align
&\sum\Sb n_2 \in {\Cal S}(y)\endSb \mu^2(n_2) \prod_{p|n_2} (c_p(z)-1) 
\Big( \frac{6}{\pi^2} x \prod_{p|n_2} \L{p}{p+1}  + O(\sqrt{x} d(n_2))\Big)
\\
= &\frac{6}{\pi^2} x {\Bbb E}(L(1,X;y)^z) + 
O\Big(\sqrt{x} \prod_{p\le y} \Big( 1+ 2|c_p(z)-1|\Big)\Big).
\\
\endalign
$$
The error term here may  be subsumed into
our estimate for the contribution of the $n_1 >1$ terms.

Thus we have an asymptotic formula for $\sumf_{|d|\le x} L(1,\chi_d;y)^z$ 
with an error 
$$
\align
&\ll x^{1-\frac{l}{8L}} y^{\frac 13} \prod_{p \le \min(4|z|+4,y)} 
\Big( 8p^{\frac 1{7L}} \Big(1-\frac {\delta}p\Big)^{-\text{Re } z} \Big) 
\prod_{\min(4|z|+4, y)\le p\le y} \Big(1+\frac{8|z|p^{1/7L}}{p}\Big)
\\
&\ll x^{1-\frac{l}{8L}} y^{\frac 13} {\Bbb E}(L(1,X;y)^{\text{Re }z}) 
\prod_{p \le \min(4|z|+4,y)} (24p^{\frac 1{7L}})
\prod_{\min(4|z|+4, y)\le p\le y} \Big(1+\frac{8|z|p^{1/7L}}{p}\Big),
\\
\endalign
$$
since ${\Bbb E}(L(1,X;y)^{\text{Re }z}) \ge \prod_{p\le \min(4|z|+4,y)} \frac 12\L{p}{p+1} (1-\delta/p)^{-\text{Re }z}$.

Suppose that $y \ge 4|z|+4$.  Then choosing the integer $l$ such 
that $2\log y > L=2^l \ge \log y$ and using the prime number 
theorem we see easily that the error above meets the 
bound prescribed in the Theorem.  
Suppose now that $y \le 4|z|+4$.   If $y \le \frac 12\log x$ then 
we choose $l=2$ and then a simple calculation gives that 
the error above is $\ll x^{1-\frac 1{40}}$, proving the Theorem in 
this case.  Lastly, when $y \le \frac 13 \log x \log_3 x$, we take 
$l = [\log_3 x]$, and again the error above
meets the bound in the Theorem.

\enddemo 

In applications, the case $y > 4|z|+4$ is the most useful case, but 
we have included the other cases mainly for the sake of 
completeness.  If the GRH is true then the bounds 
for $E(z,y)$ may be improved considerably. 
 
\proclaim{Theorem 6.2}  Let $y$, $z$ and $E(z,y)$ 
be as in Theorem 6.1 above.  Suppose that the GRH is true.  
If $y\ge 4|z|+4$ then 
$$
E(z,y) \ll x^{-\frac 1{10}} 
\exp\Big( \frac{30|z|}{\log (4|z|+4)} 
+ 10|z| \log \Big(\frac{\log y}{\log (4|z|+4)}
\Big)\Big),
$$
while if $y \le 4|z|+4$ then 
$$
E(z,y) \ll x^{-\frac{1}{10}} \exp\Big( \frac{4y}{\log y}\Big).
$$
\endproclaim 
\demo{Proof}  We follow the proof of Theorem 6.1 closely.  The 
main change is that for $n_1>1$ we have the following 
improved bound for $S(x;n_1 n_2^2)$ arising from 
Lemma 4.5:
$$
S(x;n_1n_2^2) \ll x^{\frac 34+\varepsilon} \exp(\sqrt{\log (n_1n_2^2)}) 
\ll x^{\frac 9{10}} (n_1n_2^2)^{\frac{10}{\log x}},
$$
where the last estimate follows upon distinguishing the 
cases $n_1n_2^2 \ge \exp((\log x)^2/100)$ and $n_1n_2^2 \le 
\exp((\log x)^2/100)$.  
\enddemo 

Using Theorems 6.1 and 6.2, we may now prove Theorem 3.

\demo{Proof of Theorem 3}  
Taking $A=4$ and $y= (\log x)^{180}$ in (2.1) of Proposition 2.2 
we see that for all but at most 
$x^{\frac 12}$ fundamental discriminants $|d|\le x$ 
we have
$$
L(1,\chi_d) = L(1,\chi_d;y)\Big(1 + O\Big( \frac{1}{(\log x)^{10}}\Big)\Big). 
$$
We take ${\Cal E}$ to be the 
set of exceptional discriminants $d$ for which the above 
asymptotic does not hold.  Then for a complex number $z$ with $|z| 
\le \log x$ we have
$$
\sumf\Sb |d|\le x \\ d\notin {\Cal E}\endSb 
L(1,\chi_d)^z 
= \sumf\Sb |d|\le x\\ d\notin {\Cal E}\endSb 
L(1,\chi_d;y)^z + O\Big( \frac{1}{(\log x)^9} 
\sumf\Sb |d|\le x \endSb L(1,\chi_d;y)^{\text{Re }z}\Big).
$$
In the notation of Theorem 6.1 we have that 
$$
\align
\sumf\Sb |d|\le x \\ d\notin {\Cal E}\endSb 
L(1,\chi_d)^z &= \frac{6}{\pi^2} x {\Bbb E} (L(1,X;y)^z) 
+ O\Big( x {\Bbb E}(L(1,X,y)^{\text{Re }z}) 
\Big( \frac{1+E(\text{Re }z,y)}{(\log x)^9} +E(z,y) \Big)\Big)
\\
&\hskip 1 in 
+O\Big(|{\Cal E}| \prod_{p\le y} \Big(1-\frac{\delta}{p}\Big)^{-\text{Re }z}
\Big).
\\
\endalign
$$ 
Since $|{\Cal E}| \le x^{\frac 12}$ and ${\Bbb E}(L(1,X;y)^{\text{Re }z})
\ge \prod_{p\le 4|z|+4} \L{p}{2(p+1)} (1-\delta/p)^{-\text{Re }z}$ 
we see that in the range $|z|\le 10^{-3}\log x$ 
the last error term above is 
$$
\align
&O\Big(x^{\frac 12} {\Bbb E}(L(1,X;y)^{\text{Re }z} ) \exp\Big( 
\frac{8|z|}{\log (4|z|+4)} +|z|\log \Big(\frac{\log y}{\log (4|z|+4)} 
\Big)\Big)\Big) \\
= &O(x^{\frac 23} {\Bbb E}(L(1,X;y)^{\text{Re }z})).
\\
\endalign
$$
Further, applying Theorem 6.1 we get that if $|z|\le e^{-12} \log x \log_3 x
/\log_2 x$ then $E(z,y)$ and $E(\text{Re }z,y)$ are $\ll (\log x)^{-9}$.  
If the GRH is true then Theorem 6.2 gives that 
$E(z,y)$ and $E(\text{Re }z,y)$ are $\ll (\log x)^{-9}$
in the extended range $|z|\le 10^{-3} \log x$.  
Thus we have shown that 
$$
\sumf\Sb |d|\le x\\ d\notin {\Cal E}\endSb L(1,\chi_d)^z 
= \frac 6{\pi^2} x {\Bbb E}(L(1,X;y)^z) 
+ O\Big( x \frac{{\Bbb E}(L(1,X;y)^{\text{Re }z})}{(\log x)^9}\Big),
$$ 
in the appropriate conditional and unconditional ranges for $z$.  
Observe that ${\Bbb E}(L(1,X)^z) = {\Bbb E}(L(1,X;y)^z) 
\prod_{p>y} (1+O(|z|^2/p^2)) = {\Bbb E}(L(1,X;y)^{z}) (1+O((\log x)^{-178}))$, 
and further that 
${\Bbb E}(L(1,X;y)^{\text{Re }z })\le {\Bbb E}(L(1,X)^{\text{Re }z})
$, and so Theorem 3 follows.

\enddemo

\head 7. The distribution function: Proof of Theorem 1 \endhead

\noindent As with Proposition 1 we prove only the 
statements involving $\Phi_{\Bbb R}$, the corresponding results 
for $\Psi_{\Bbb R}$ are established similarly. 
Let $\lambda >0$ be a real number, and let $N\ge 1$ be 
an integer.  Observe that for any $y>0$ and $c>0$ we have 
$$
\frac{1}{2\pi i} \int_{c-i\infty}^{c+i\infty}  y^s \Big(\frac{e^{\lam s}-1}
{\lam s}\Big)^N \frac {ds}{s} 
= \frac{1}{\lam^N } \int_0^{\lam}\ldots \int_0^\lam 
\frac{1}{2\pi i} \int_{c-i\infty}^{c+i\infty} 
(ye^{t_1+\ldots+t_N})^s \frac{ds}{s} dt_1 \ldots dt_N 
$$
so that by Perron's formula we get 
$$
\frac{1}{2\pi i} \int_{c-i\infty}^{c+i\infty}  y^s \Big(\frac{e^{\lam s}-1}
{\lam s}\Big)^N \frac {ds}{s} 
\cases
=1 &\text{if  } y\ge 1 \\
\in [0,1] &\text{if  } 1\ge y\ge e^{-\lam N}\\
=0 &\text{if  } e^{-\lam N} \ge y. \\
\endcases
$$

Write $s=k+it$ from now on, where $k$ is positive.  
We put 
$$
I_1 = \frac{1}{2\pi i} \int_{(k)} \frac{\pi^2}{6 x} 
\sumf\Sb |d|\le x\\ d\notin {\Cal E}\endSb 
L(1,\chi_d)^s (e^{\gamma} \tau)^{-s} \Big(\frac{e^{\lam s}-1}{\lam s}\Big)^{N}
\frac{ds}{s}, 
$$
where ${\Cal E}$ is the set in Theorem 3, so that $|{\Cal E}|\le 
x^{\frac 12}$.  
>From the above we know that 
$$
{\Phi_x}(\tau) +O(x^{-\frac 12}) \le I_1
\le {\Phi_x}( \tau e^{-\lam N}) +O(x^{-\frac 12}).   
\tag{7.1} 
$$
Further if we set 
$$
I_2 = \frac{1}{2\pi i} \int_{(k)} {\Bbb E}(L(1,X)^s)
(e^{\gamma} \tau)^{-s} \Big(\frac{e^{\lam s}-1}{\lam s}\Big)^{N}
\frac{ds}{s}, 
$$ 
then we have that 
$$
\Phi(\tau) \le I_2 \le \Phi(\tau e^{-\lam N} ) . 
\tag{7.2}
$$

We shall now estimate $|I_1-I_2|$.  We put $T:= 10^{-6} \log x 
\log_3 x/\log_2 x$ if we are to prove the unconditional parts of 
Theorem 1, and $T:= 10^{-4} \log x$ if we are to 
prove the GRH part of Theorem 1.  We take $k=k_\tau$ as in Theorem 
3.1 (thinking of $y$ there as $\to \infty$; thus $I(k)=\gamma+\log \tau$), 
so that when $\tau \le \log T + C_1 -20$ we have $k_\tau 
= e^{\tau -C_1} \{ 1+O(1/\tau)\} \le e^{-15}T$.  Observe that 
$$
|I_1 -I_2| \ll (e^{\gamma} \tau)^{-k} 
\Big( \int_{|t|\le T} +\int_{|t|> T}\Big) 
\Big|\frac{\pi^2}{6 x} \sumf\Sb |d|\le x\\ d\notin {\Cal E}
\endSb L(1,\chi_d)^s - {\Bbb E}(L(1,X)^s) \Big| \Big|\frac{e^{\lam s}-1}
{\lam s}\Big|^N \Big|\frac{ds}{s}\Big|.
\tag{7.3}
$$

From Theorem 3.1 we know that $\Phi(\tau) 
\asymp \left( \sqrt{(\log k)/k} \right) \left( {\Bbb E}(L(1,X)^{k})/(e^{\gamma} \tau)^k \right)$, and 
so applying Theorem 3 we get that when $|t|\le T$ 
$$
(e^{\gamma} \tau)^{-k} 
\Big|\frac{\pi^2}{6 x} \sumf\Sb |d|\le x\\ d\notin {\Cal E}
\endSb L(1,\chi_d)^s - {\Bbb E}(L(1,X)^s) \Big| 
\ll \frac{1}{(\log x)^9} \frac{{\Bbb E} (L(1,X)^{k})}{(e^{\gamma}\tau)^k} 
\ll \frac{\Phi(\tau) }{(\log x)^8}.
$$
Further note that $|(e^{\lam s}-1)/(\lam s)| \le 1+ e^{\lam k}$, 
which is easily seen by looking at the cases $|\lam s|\le 1$ and 
$|\lam s|>1$.  Hence it follows that the $|t|<T$ integral 
contributes to (7.3) an amount 
$$
\ll \frac{\Phi(\tau)}{(\log x)^8} (1+e^{\lam k})^N \log T 
\ll (1+e^{\lam k})^N \frac{\Phi(\tau)}{(\log x)^7}.
$$

Next note that, by Theorem 3,
$$
\align
(e^{\gamma} \tau)^{-k} 
\Big|\frac{\pi^2}{6 x} \sumf\Sb |d|\le x\\ d\notin {\Cal E}
\endSb L(1,\chi_d)^s - {\Bbb E}(L(1,X)^s) \Big| 
&\ll (e^{\gamma} \tau)^{-k} 
\Big(\frac{\pi^2}{6 x} \sumf\Sb |d|\le x\\ d\notin {\Cal E}
\endSb L(1,\chi_d)^k +   
 {\Bbb E}(L(1,X)^k)\Big) 
\\
&\ll \frac{{\Bbb E}(L(1,X)^k)}{(e^{\gamma} \tau)^k} \ll  
\Phi(\tau)  \sqrt{\log x} . 
\\
\endalign
$$
Hence the contribution to (7.3) from the $|t|>T$ integral is 
$$
\ll \sqrt{\log x} 
\Phi(\tau) \int_{|t|>T} \frac{(1+e^{\lam k})^N}{(\lam |t|)^N} 
\frac{dt}{|t|} 
\ll  \sqrt{\log x} 
\Phi(\tau) \Big(\frac{1+e^{\lam k}}{\lam T} \Big)^N.
$$

We now choose $\lam = 4e^{10}/T$, and $N=[\log_2 x]$.
>From the above estimates we then conclude that 
$$
|I_1 -I_2| \ll {\Phi(\tau)} (1+e^{\lam k})^N 
\Big( \frac{1}{(\log x)^7} + \frac{\sqrt{\log x}}{(\lam T)^N}\Big) 
\ll \frac{\Phi(\tau)}{(\log x)^5}.
$$
Using this in (7.1), and (7.2) we deduce that  
$$
{\Phi_x}(\tau)  \le \Phi(\tau e^{-\lam N}) + O\Big(x^{-\frac 12} + 
  \frac{\Phi(\tau)}{(\log x)^5}\Big),  
$$ 
and that 
$$
{\Phi_x}( \tau e^{-\lam N}) \ge  \Phi(\tau) + O\Big( x^{-\frac 12} + 
  \frac{\Phi(\tau)} {(\log x)^5}   \Big).
$$
In the lower bound above, replace $\tau$ by $\tau e^{\lam N}$.  Then 
we get that 
$$
\Phi(\tau e^{\lam N})  (1+O((\log x)^{-5}) ) 
\le {\Phi_x}(\tau) + O( x^{-\frac 12})
\le \Phi(\tau e^{-\lam N}) ( 1+O((\log x)^{-5}) ) ,
$$
since $\Phi(\tau) \leq \Phi(\tau e^{-\lam N})$.
We now observe that $\Phi(\tau e^{\pm\lam N})=\Phi(\tau) 
\exp(O(\lam N e^{\tau}))$ which follows directly from Proposition 1 if 
$\lam N \le e^{-\tau}$,  and also in the range $\lam N \le 1/\tau$ upon 
iterating the estimate in Proposition 1.  
Theorem 1 now follows upon recalling our 
choices of $\lam$ and $N$.

\head 8. The distribution function: Proof of Theorem 4\endhead 

\noindent To prove Theorem 4 we shall first show that for most 
discriminants $d$ we may approximate $L(1,\chi_d)$ by a  short Euler 
product.  Then we shall use Theorems 6.1 and  6.2
to compare the distribution 
of these short Euler products with the random Euler products of 
Section 3.

\proclaim{Lemma 8.1}  
Let $y^{100} \ge z\ge y$ be real numbers, with $y$ large.  
Then for integers $k$, $l\ge 2$ 
$$
\sumf_{|d|\le x} \left| \sum_{y\le p\le z} \frac{\L{d}{p}}{p} \right|^{2k} 
\ll x \Big(\frac{2k}{y\log y}\Big)^{k} 
+ x^{1-\frac{l}{8L}} y^{35+30\frac{k}{L}} e^{6k}
$$
uniformly, where $L=2^l$.  
Further if the GRH is assumed then 
$$
\sumf_{|d|\le x} \left| \sum_{y\le p\le z} \frac{\L{d}{p}}{p} \right|^{2k} 
\ll x \Big(\frac{2k}{y\log y}\Big)^{k} + x^{\frac 45} 
\exp(\sqrt{2k\log z}) \Big(\frac{\log (ez/y)}{\log y}\Big)^{2k}.
$$
\endproclaim

\demo{Proof}  The quantity we seek to estimate is 
$$
\sum\Sb p_1, \ldots, p_{2k}\\ y\le p_i \le z\endSb 
\frac{1}{p_1\cdots p_{2k}} \sumf_{|d|\le x} \L{d}{p_1\cdots p_{2k}}.
$$
The terms where $p_1\cdots p_{2k} =\square$ (so that the $p_i$ occur in pairs)
contribute 
$$
\ll x \frac{(2k)!}{2^k k!} \Big(\sum_{y\le p \le z} \frac{1}{p^2} 
\Big)^k \ll x \Big(\frac{2k}{y\log y} \Big)^k,
$$
using Stirling's formula and $\sum_{p>y} \frac{1}{p^2} \le \frac{2}{y\log y}$.
By Lemma 4.2, the terms where $p_1\ldots p_{2k} \neq \square$ 
contribute 
$$
\align
&\ll \sum\Sb p_1, \ldots, p_{2k}\\ y\le p_i \le z\endSb 
\frac{1}{p_1\cdots p_{2k}} 
x^{1-\frac{l}{8L}}  y^{35} (y^{200k})^{\frac{1}{7L}} 
(2^{2k})^{\frac{l^2}{L}} \prod_{i=1}^{2k} 
\Big(1+\frac{1}{p_i^{1-l/8L}}\Big)
\\
&\ll x^{1-\frac{l}{8L}} y^{35+30\frac{k}{L}} \Big(\sum_{y<p<y^{100}} 
\frac{e}{p}\Big)^{2k} 
\ll 
 x^{1-\frac{l}{8L}} y^{35+30\frac{k}{L}} e^{6k}, 
\\
\endalign
$$
which gives the first part of the Lemma.  If the 
GRH is assumed then we may estimate the contribution of the 
terms $p_1\ldots p_{2k} \neq \square$ using Lemma 4.5.  
This gives for these terms an amount
$$
\ll 
x^{\frac 45} \exp(\sqrt{2k\log z}) 
\Big( \sum_{y\le p\le z} \frac 1p\Big)^{2k} \ll 
x^{\frac 45} \exp(\sqrt{2k\log z}) \Big(\frac{\log (ez/y)}{\log y}\Big)^{2k},
$$ 
proving the Lemma.  
\enddemo 

\proclaim{Proposition 8.2}  Let $(\log_2 x)^{10} \ge A\ge e$ be 
a real number, put $y_0= A^2 \log x \log_3 x$, and $y_1 =eA^3 \log x \log_2 x$.
Then there is a positive constant $c$ such that 
$$
L(1,\chi_d) = L(1,\chi_d;y_0) \Big( 1+ O\Big(\frac{1}{A\log_2 x} +
\frac{1}{\log x}\Big)\Big),
$$
for all but $x \exp(-c\log x \log_3 x/ \log_2 x)$ fundamental 
discriminants $|d| \le x$.  If the GRH is true then 
$$
L(1,\chi_d) = L(1,\chi_d;y_1) \Big( 1+ O\Big(\frac{1}{A\log_2 x} +
\frac{1}{(\log_2 x)^2}\Big)\Big),
$$
for all but $x^{\frac{24}{25}}$ fundamental discriminants $|d|\le x$.  
\endproclaim

\demo{Proof} Setting $z=(\log x)^{80}$, we know that 
$L(1,\chi_d) = L(1,\chi_d;z)
(1+O(1/\log x))$ for all but at most $x^{\frac 15}$ discriminants $|d|\le x$
(by (2.1) with $A=10$). 

Since $L(1,\chi_d;z) =L(1,\chi_d;y_0) 
\exp(\sum_{y_0\le p\le z} (\chi_d(p)/p+O(1/p^2)))$, we 
see that the first part of the Proposition follows from showing that 
$$
\Big|\sum_{y\le p\le z} \frac{\chi_d(p)}{p} \Big| \ge \frac{1}{A\log y} 
\ \ \text{for } \ll x\exp\Big(-c \frac{\log x\log_3 x}{ \log_2 x}\Big) $$
fundamental discriminants $|d|\le x$. 
Choosing $k=[\log x \log_3 x/(333 \log_2 x)]$ 
and $l=[\log_3 x/(2\log 2)]$ in Lemma 8.1, 
and keeping in mind $y\ge A^2 \log x \log_3 x$, 
we deduce that 
$$
\sumf_{|d|\le x} \Big| \sum_{y\le p\le z} 
\frac{\chi_d(p)}{p}\Big|^{2k } 
\ll x \Big(\frac{1}{10 A  \log y}\Big)^{2k} 
+ x^{1-\frac{1}{\sqrt{\log_2 x}}}.
$$
>From this the first part of the Proposition follows easily.

To prove the second part, as above it suffices to show that 
$$
\Big|\sum_{y_1\le p\le z} \frac{\L{d}{p}}{p} \Big|
 \le \frac{4}{A\log_2 x}, \qquad
\text{for all but }\le x^{\frac {24}{25}} 
\text{ exceptions } |d|\le x.
$$
Put $y_r=e^{r-1} y_1$, and $z_r =e^{r} y_1$ where $r=1$, $\ldots$, 
$R$ where $R:=[80\log_2 x -\log y_1]$, and set 
$y_{R+1}=e^{R} y$, and $z_{R+1}=z$.  
For each $1\le r \le R+1$, we use the conditional part of 
Lemma 8.1 appropriately, choosing 
the exponent $k_r = [\log x/(50 \log(A e^{r/3+1}))]$.  
Keeping in mind that $y_r = A^3 e^r \log x \log_2 x$, and that 
$z_r \le \log^{80} x$, we deduce easily that 
$$
\Big|\sum_{y_r \le p\le z_r } \frac{\L dp}{p} \Big| 
\le \frac{1}{e^{r/3} A \log_2 x }, 
\qquad \text{ for all but }\ll x^{\frac{19}{20}+o(1)} 
\text{ exceptions } |D|\le x.
$$
Summing this over $1\le r\le R+1$, we obtain the second part 
of the Proposition. 

\enddemo

Keep $y_0$ and $y_1$ as in Proposition 8.2.  Then applying Theorem 6.1
we see that uniformly for $|z| \le \log x \log_3 x/(10^{6} \log A)$ 
$$
\frac{\pi^2}{6x} \sumf_{|d|\le x} L(1,\chi_d;y_0)^z 
= {\Bbb E}(L(1,X;y_0)^z) + O\Big( \frac{{\Bbb E}(L(1,X;y_0)^{\text{Re }z} ) }
{(\log x)^{10}}\Big).  \tag{8.1} 
$$
Thus upon using the argument of Section 7 together with 
Corollary 3.3 we easily obtain 
that uniformly in $\tau \le \log_2 x +\log_4 x -\log_2 A -20$ we 
have 
$$
\frac{\pi^2}{6x} \sumf\Sb |d|\le x \\ L(1,\chi_d;y_0)\ge e^{\gamma}\tau\endSb 
1 = \exp\Big( -\frac{e^{\tau-C_1}}{\tau} \Big( 1+O\Big(\frac{1}{\tau} + 
\frac{e^{\tau}}{y_0}\Big)\Big)\Big).
$$
Combining this with Proposition 8.2 above, we get the unconditional part 
of Theorem 4 for the frequency of large values.  The frequency of small 
values is obtained similarly.  As for the part of Theorem 4 conditional on GRH 
we may use the above argument, applying Theorem 6.2 instead of Theorem 
6.1.  Then (8.1) holds, with $y_0$ replaced by $y_1$, in the 
wider range $|z|\le \log x \log_2 x/(10^6 \log A)$, and hence we get the 
appropriately stronger conclusion.

\head 9. Extreme values of $L(1,\chi_d)$ under GRH: Proof of Theorem 5a
\endhead

\proclaim{Proposition 9.1}  Suppose the GRH is true.  
Let $z\ge 2$ be a real number, and let 
$P(z)=\prod_{p\le z} p = e^{z+o(z)}$.  Let $\epsilon_p = \pm 1$ for each 
prime $p\le z$, and let ${\Cal P}(x,\{\epsilon_p\})$ 
denote the set of primes $q\le x$ such that $\L{p}{q} =\epsilon_p$ for 
each $p\le z$.  Then, for $z\ll x^{1/2}$,
$$
\sum_{q\in {\Cal P}(x,\{\epsilon_p\})} \log q 
= \frac{x}{2^{\pi(z)}} + O(x^{\frac 12} \log^2 (xP(z))), \tag{9.1}
$$
and 
$$
\sum_{q\in {\Cal P}(x,\{\epsilon_p\})} L(1,(\tfrac{\cdot}{q})) \log q 
= \frac{x\zeta(2)}{2^{\pi(z)}} \prod_{p\le z} \Big(1+\frac{\epsilon_p}{p}
\Big) 
+ O(x^{\frac 12 }\log^3 (xP(z))). \tag{9.2}
$$
\endproclaim

\demo{Proof}  We shall only demonstrate (9.2): (9.1) 
is similar and simpler.
For $l|P(z)$ put $\epsilon_l=\prod_{p|l} \epsilon_p$.  Note that for $q\le x$
$$
\frac{1}{2^{\pi(z)}} \sum_{l|P(z)} \epsilon_l \L{l}{q} = 
\cases 
1 &\text{if} \qquad q\in {\Cal P}(x,\{\epsilon_p\})\\
0 &\text{if otherwise}.\\
\endcases
$$
By partial summation we easily see that 
$$
L(1,(\tfrac{\cdot}{q})) = \sum_{n\le N} \frac{(\tfrac{n}{q})}{n}
+O\Big(\frac{q}{N}\Big). \tag{9.3}
$$
Using (9.3) with $N=x^2 P(z)^2$, it follows that the LHS of (9.2) equals
$$
\frac{1}{2^{\pi(z)}} \sum_{l|P(z)} \epsilon_l \sum_{n\le x^2 P(z)^2} 
\frac 1n \sum_{q\le x} \L{ln}{q} \log q + O(1/P(z)^2). \tag{9.4}
$$

If $ln=\square$ then the inner sum over $q$ above is 
$\psi(x) + O(\sqrt{x} +\sum_{p|ln} \log p) = x+O(x^{\frac 12}\log^2 x +\log (xP(z)))
=x+O(x^{\frac 12}\log^2 x)$,
since RH is assumed.  Moreover we may write $n=lm^2$ since $l$ is squarefree.  
Thus the contribution of these terms to (9.4) is
$$
\align
&\frac{1}{2^{\pi(z)}} \sum_{l|P(z)}{\epsilon_l} 
\sum_{m\le xP(z)/\sqrt{l}} \frac{1}{lm^2} 
(x+O(x^{\frac 12}\log^2 x))
\\
=& \frac{\zeta(2)}
{2^{\pi(z)}} \prod_{p\le z} \Big(1+\frac{\epsilon_p}{p}\Big) 
(x+O(x^{\frac 12}\log^2 x))
\\
\endalign
$$
This accounts for the main term in (9.2), with an acceptable error term.

If $ln\neq \square$ then $\L{ln}{\cdot}$ is a non-principal character $\chi$ 
say, of conductor $ln$ or $4ln$.  Thus the inner sum over $q$ in (9.4) 
is now $\psi(x,\chi) +O(x^{\frac 12}) \ll x^{\frac 12} \log^2 (xP(z))$,
by GRH.  It follows that the contribution of these terms to (9.4) 
is 
$$
\ll  x^{\frac 12} \log^2 (xP(z)) \frac{1}{2^{\pi(z)}} \sum_{l|P(z)} 
\sum_{n\le x^2 P(z)^2} \frac 1n \ll x^{\frac 12} \log^3 (xP(z)).
$$
This proves the Proposition.
\enddemo

We are now ready to prove Theorem 5a. 
Let $z$ be such that $2^{\pi(z)} \le x^{\frac 12-3\epsilon}$.  
Then for any choice of $\epsilon_p =\pm 1$ for $p\le z$, we get by Proposition
9.1 that
$$
\sum_{q\in {\Cal P}(x,\{\epsilon_p\})} \log q 
= \frac{x}{2^{\pi(z)}} (1+O(x^{-2\epsilon})),
$$
and 
$$
\sum_{q\in {\Cal P}(x,\{\epsilon_p\})} L(1,(\tfrac{\cdot}{q})) \log q 
= \frac{x\zeta(2)}{2^{\pi(z)}} \prod_{p\le z} \Big(1+\frac{\epsilon_p}{p}
\Big) (1+ O(x^{-2\epsilon})). 
$$
Since $L(1,\L{\cdot}{q}) \le 2\log q$ for large $q$ (take $N=q$ in (9.3)),
we deduce easily from the above estimates that for $\delta>0$ 
there are 
$$
\ge \frac{\delta}{2\log^2 x} \frac{x}{2^{\pi(z)}} + 
O\Big(\frac{x^{1-\epsilon}}{2^{\pi(z)}}\Big)
\tag{9.5}
$$
primes $q\in {\Cal P}(x,\{\epsilon_p\})$ with 
$$
L(1,(\tfrac{\cdot}{q})) \ge \zeta(2) \prod_{p\le z} \Big(1+
\frac{\epsilon_p}{p}\Big) -\delta.
$$
Similarly we deduce that the sum of $\log q$, over those 
primes $q\in {\Cal P}(x,\{\epsilon_p\})$ with 
$$
L(1,(\tfrac{\cdot}{q})) \le \zeta(2) \prod_{p\le z} \Big(1+
\frac{\epsilon_p}{p}\Big) (1+\delta),
$$
is 
$$
\ge \frac{\delta }{(1+\delta)} 
\frac{x}{2^{\pi(z)}} 
+ O\Big(\frac{x^{1-\epsilon}}{2^{\pi(z)}}\Big). \tag{9.6}
$$

Take $z=\log x \log \log x /(2\log 2+ 10\epsilon)$,  $\delta=\epsilon$,
and $\epsilon_p=1$ for all $p\le z$.  By (9.5) that there are 
$\gg x^{\frac 12}$ 
primes $q\le x$ such that 
$$
L(1,(\tfrac{\cdot}{q})) \ge \zeta(2) \prod_{p\le z} \Big(1+
\frac{1}{p}\Big) -\epsilon 
\ge e^{\gamma} (\log \log x+\log \log \log x -\log(2\log 2) -100\epsilon),
$$
by using the prime number theorem, which implies the first part 
of Theorem 5a. The second part 
follows from the analogous argument using (9.6) with  
the same $z$, and with $\delta= \epsilon/\log x$, and $\epsilon_p=-1$ for all 
$p\le z$.

\head 10.  Large values of $L(1,\chi_d)$: Proof of Theorem 5b \endhead

\noindent Let $x$ be large, and put $L= [100\log_3 x]$, and 
$z= \log x \log_2 x/(10L)$.  

\proclaim{Lemma 10.1} There are at least $x^{\frac{1}{3L}}$ pairwise 
coprime integers of the form $pq$ with $p$ and $q$ primes below 
$x^{\frac{1}{2L}}$ such that $\L{\ell}{pq} =1$ for all primes $\ell \le z$.
\endproclaim 

\demo{Proof}  For each prime $\log^2 x \le p \le x^{\frac 1{2L}}$ let 
$\delta(p)$ denote the vector $(\L{2}{p}, \ \L{3}{p}, \ \ldots, \ \L{
\ell_{\text{max}}}{p})$ where $\ell_{\text{max}}$ denotes the largest 
prime below $z$.  Note that the total number of possible vectors 
is $2^{\pi(z)}$. Also observe that if $p$ and $q$ have the same 
vector $\delta(p)=\delta(q)$ then $\L{\ell}{pq}=1$ for all primes 
$\ell \le z$.  

Given a vector $\delta$ of $\pm 1$'s let $N(\delta)$
denote the number of primes $p$ with $\delta(p)=\delta$.  
By considering the products of these primes taken two at a time we 
get $\ge [N(\delta)/2]$ coprime integers of the form $pq$, satisfying 
$\L{\ell}{pq}=1$ for all $\ell \le z$.  Thus the total number 
of such $pq$ exceeds 
$$
\sum_{\delta} (N(\delta)/2-1) = \frac 12 (\pi(x^{\frac{1}{2L}}) -\pi(\log^2 x) )
- 2^{\pi(z)} \ge x^{\frac{1}{3L}},
$$
proving the Lemma.
\enddemo

\proclaim{Lemma 10.2} Let $d_1$, $\ldots$, $d_L$ be any $L$ numbers 
constructed in Lemma 10.1.  Then there exists a square-free integer $d\le x$ 
which is the product of at least $L/3$ of these $d_i$'s such that 
$$
\L{\ell}{d} = 1, \qquad \text{for all primes } \ell \le z; \qquad 
\text{and} \ \sum_{z\le p \le \exp((\log x)^{3})} \frac{\L{p}{d}}{p } 
\ge -\frac{1}{\log_2 x}. 
$$
\endproclaim 

\demo{Proof}  If $d$ is a product of distinct $d_i$ let $\nu(d)$ 
denote the number of $d_i$'s involved in this product.  Note that 
$$
\sum_{z\le p \le \exp((\log x)^{3})} \frac{1}{p} \prod_{i=1}^{L } 
\Big( 1+\L{p}{d_i}\Big) \ge 0,
$$
and so 
$$
\sum\Sb d\\ \nu(d) \ge L/3\endSb \sum_{z\le p\le \exp((\log x)^{3})} 
\frac{1}{p} \L{p}{d} 
\ge - \sum\Sb d\\ \nu(d)\le L/3 \endSb  \sum_{z\le p\le \exp((\log x)^{3})} 
\frac{1}{p} 
\ge -10 \log_2 x \sum\Sb d\\ \nu(d)\le L/3 \endSb 1.
$$
Now $\# \{ d:\ \nu(d)\ge L/3\} \leq \#\{ d\} = 2^L$, and 
$\# \{ d:\ \nu(d)\le L/3\} =\sum_{j\le L/3} \binom{L}{j} 
\le (1.9)^L$ using Stirling's formula.  It follows that 
there exists $d$ with $\nu(d)\ge L/3$ and 
$$
\sum_{z\le p\le \exp((\log x)^{100})} 
\frac{1}{p} \L{p}{d}  \ge - 10 (\log_2 x) (0.95)^L \ge -\frac{1}{\log_2 x},
$$
as desired.  

\enddemo

\proclaim{Lemma 10.3}  There are $\ge x^{\frac{2}{19}}$ integers 
$d\le x$ as in Lemma 10.2. 
\endproclaim

\demo{Proof}  Counted with multiplicity the number of $d$'s 
constructed in Lemma 10.2 is $\binom{[x^{\frac{1}{3L}}]}{L}$.  The number 
$d$ is counted $\binom{[x^{\frac{1}{3L}}]-\nu(d)}{L-\nu(d)} \geq \binom{[x^{\frac{1}{3L}}]-[L/3]}{L-[L/3]}$ 
times, since each such $d$ has $\nu(d)\leq L/3$. So we are left with at least 
$\binom{[x^{\frac{1}{3L}}]}{L}/\binom{[x^{\frac{1}{3L}}]-[L/3]}{L-[L/3]}$
 distinct $d$, which suffices.
\enddemo

\demo{Proof of Theorem 5b}  If $d$ is such that there is no Landau-Siegel 
character $\pmod {d}$ then by a simple application of 
the prime number theorem for arithmetic progressions 
we see that 
$$
\log L(1,(\tfrac{\cdot}{d})) 
=\sum\Sb 2\le n\le \exp(\log^3 x) \endSb 
\frac{\Lambda(n)}{n\log n} \L{n}{d} 
+ O\Big( \frac{1}{\log x}\Big), 
$$
so that 
$$
L(1,(\tfrac{\cdot}{d})) 
= \prod_{p\le \exp((\log x)^{3})} 
\Big(1-\frac{\L{p}{d}}{p}\Big)^{-1} \Big(1+O\Big(\frac{1}{\log x}\Big)\Big).
$$
Use this for the $d$ in Lemma 10.3 which are not Landau-Siegel moduli, 
and the result follows.

\enddemo

\enddocument